\documentclass[oupthm]{my}

\usepackage{graphicx}
\usepackage{stackengine}
\usepackage{scalerel}
\usepackage{tensor}       %\tensor[_n]{\braket{j|\psi}}{_{1\ldots n}}
\usepackage{array}
\usepackage{makecell}
\newcolumntype{x}[1]{>{\centering\arraybackslash}p{#1}}
\usepackage{tikz}
\usepackage{pgfplots}
\stackMath

\makeatletter

\tikzset{meter/.append style={draw, inner sep=10, rectangle, font=\vphantom{A}, minimum width=30, line width=.8, path picture={\draw[black] ([shift={(.1,.3)}]path picture bounding box.south west) to[bend left=50] ([shift={(-.1,.3)}]path picture bounding box.south east);\draw[black,-latex] ([shift={(0,.1)}]path picture bounding box.south) -- ([shift={(.3,-.1)}]path picture bounding box.north);}}}
\definecolor{Blues5seq1}{RGB}{239,243,255}
\definecolor{Blues5seq2}{RGB}{189,215,231}
\definecolor{Blues5seq3}{RGB}{107,174,214}
\definecolor{Blues5seq4}{RGB}{49,130,189}
\definecolor{Blues5seq5}{RGB}{8,81,156}

\definecolor{Greens5seq1}{RGB}{237,248,233}
\definecolor{Greens5seq2}{RGB}{186,228,179}
\definecolor{Greens5seq3}{RGB}{116,196,118}
\definecolor{Greens5seq4}{RGB}{49,163,84}
\definecolor{Greens5seq5}{RGB}{0,109,44}

\definecolor{Reds5seq1}{RGB}{254,229,217}
\definecolor{Reds5seq2}{RGB}{252,174,145}
\definecolor{Reds5seq3}{RGB}{251,106,74}
\definecolor{Reds5seq4}{RGB}{222,45,38}
\definecolor{Reds5seq5}{RGB}{165,15,21}

\usepackage[most,breakable]{tcolorbox}
\makeatletter
\renewenvironment{boxed}{\begingroup\@ifnextchar\bgroup\boxed@gobegin\boxed@gobegin@empty}{\end{tcolorbox}\endgroup}
\def\boxed@gobegin#1{\def\@tempa{#1}\def\@tempb{orange}\ifx\@tempa\@tempb\begin{tcolorbox}[colback=red!15,colframe=orange!70,breakable,enhanced]\else\begin{tcolorbox}[colback=Blues5seq1,colframe=Blues5seq5,breakable,enhanced]\fi}
\def\boxed@gobegin@empty{\begin{tcolorbox}[colback=Blues5seq1,colframe=Blues5seq5,breakable,enhanced]}
\makeatother
\newcommand{\ba}[1]{\begin{align}#1\end{align}}

\newcommand{\Par}[1]{\left( #1 \right)}
\newcommand{\Brac}[1]{\left[ #1 \right]}

\newcommand{\Sp}{\operatorname{Sp}}
\newcommand{\SpSm}{\operatorname{SpS}}

\newcommand{\OrSp}{\operatorname{OSp}}
\newcommand{\Or}{\operatorname{O}}

\newcommand{\ran}{\operatorname{range}}
\newcommand{\Ker}{\operatorname{ker}}

\newcommand{\Pd}{\operatorname{Pd}}
\newcommand{\Psd}{\operatorname{Psd}}
\newcommand{\diag}{\operatorname{diag}}
\newcommand{\Span}{\operatorname{span}}
\newcommand{\rank}{\operatorname{rank}}

\usepackage{xcolor}
\usepackage[normalem]{ulem}
\usepackage{dsfont}
\usepackage{cite}
\usepackage[numbers, sort&compress]{natbib}
\usepackage{footnote}
\usepackage{fixfoot}
\usepackage{dsfont}
\usepackage{bigints}
\usepackage{tikz}
\usepackage{geometry}
\usepackage{enumitem}   
\newgeometry{bottom=20mm} 

\usepackage{float}
\usepackage{graphicx}
\usepackage{multicol}
\usepackage{amsmath,amssymb,amsfonts, mathtools}
\usepackage{mathrsfs}

\usepackage{rotating}
\usepackage[title, toc]{appendix}
\usepackage[numbers]{natbib}
\usepackage{changepage,afterpage}
\usepackage{etruscan}
\usepackage{tikz, siunitx, newunicodechar}
\usepackage{booktabs,multirow}
\tcbuselibrary{skins,breakable}
\usepackage{xcolor}
\usepackage{lipsum}
\usepackage{longtable}

\usepackage{bm,bbm}
\usepackage{braket}

\usepackage[utf8]{inputenc}
\usepackage{array}
\usepackage{colortbl}

\allowdisplaybreaks[4]

\usepackage[T3,T1]{fontenc}
\DeclareSymbolFont{tipa}{T3}{cmr}{m}{n}
\DeclareMathAccent{\invbreve}{\mathalpha}{tipa}{16}

\theoremstyle{oupplain}
\newtheorem{theorem}{Theorem}[section]
\newtheorem{lemma}[theorem]{Lemma}
\newtheorem{proposition}[theorem]{Proposition}
\newtheorem{corollary}[theorem]{Corollary}
\theoremstyle{oupdefinition}
\newtheorem{definition}{Definition}[section]
\theoremstyle{oupremark}
\newtheorem{remark}[theorem]{Remark}

\theoremstyle{oupproof}

\definecolor{ashgrey}{rgb}{0.7, 0.75, 0.71}

\newenvironment{proof}{
	\tcolorbox[blanker,breakable,left=5mm,parbox=false,
    before upper={\parindent15pt},
    after skip=10pt,
	borderline west={0.7mm}{0pt}{Blues5seq1}]
    {\noindent{\it \textbf{Proof:}}}
}{
    \textcolor{black}{\hbox{}\nobreak\hfill$\blacksquare$} 
    \endtcolorbox
}

\numberwithin{equation}{section}

\articletype{RESEARCH ARTICLE}

\begin{document}
\setlength{\abovedisplayskip}{12pt}
\setlength{\belowdisplayskip}{12pt}

\begin{Frontmatter}
\title{On generalization of Williamson's theorem to real symmetric matrices}

    \author{Hemant K. Mishra \hspace{-0.10cm} \thanks{School of Electrical and Computer Engineering, Cornell University, Ithaca, New York~14850, USA} \thanks{Department of Mathematics and Computing, Indian Institute of Technology (ISM) Dhanbad, Jharkhand 826004, India;\email{hemantmishra1124@iitism.ac.in}} }

    \abstract 
        {   
           Williamson's theorem states that if $A$ is a $2n \times 2n$ real symmetric positive definite matrix then there exists a $2n \times 2n$ real symplectic matrix $M$ such that $M^{\top} A M=D \oplus D$, where $D$ is an $n \times n$ diagonal matrix with positive diagonal entries known as the symplectic eigenvalues of $A$.
        The theorem is known to be generalized to $2n \times 2n$ real symmetric positive semidefinite matrices whose kernels are symplectic subspaces of $\mathds{R}^{2n}$, in which case, some of the diagonal entries of $D$ are allowed to be zero.
        In this paper, we further generalize Williamson's theorem to $2n \times 2n$ real symmetric matrices by allowing the diagonal elements of $D$ to be any real numbers, and thus extending the notion of symplectic eigenvalues to real symmetric matrices.
        Also, we provide an explicit description of symplectic eigenvalues, construct symplectic matrices achieving Williamson's theorem type decomposition, and establish perturbation bounds on symplectic eigenvalues for a class of $2n \times 2n$ real symmetric matrices denoted by $\operatorname{EigSpSm}(2n)$.
        % The set $\operatorname{EigSpSm}(2n)$ contains the set of $2n \times 2n$ real symmetric positive semidefinite matrices whose kernels are symplectic subspaces of $\mathds{R}^{2n}$.
        Our perturbation bounds on symplectic eigenvalues for $\operatorname{EigSpSm}(2n)$ generalize known perturbation bounds on symplectic eigenvalues of positive definite matrices given by Bhatia and Jain \textit{[J. Math. Phys. 56, 112201 (2015)]}.
        }

\end{Frontmatter}

\begin{boxed}
    \tableofcontents
\end{boxed}

\section{Introduction}
    Williamson's theorem contains germs of modern developments in symplectic topology.
    It facilitates an immediate proof of Gromov's non-squeezing theorem in the linear case \cite{gromov1985pseudo}, which is one of the most important theorems in symplectic geometry.
    Also known as Williamson's decomposition, the theorem is fundamental in developing the theory of bosonic Gaussian states in quantum information \cite{serafini2003symplectic, pereira2021symplectic, nicacio2021williamson, vsafranek2015quantum}.
    In the recent years, Williamson's theorem has attracted much attention of mathematicians and physicists,
    and it has become a topic of intense study in matrix analysis \cite{bhatia2015symplectic, HIAI2018129, mishra2020first, bhatia2020schur, bhatia_jain_2021, jain2021sums, jm, paradan2022horn, mishra2023, sags_2021, huang2023, son2022symplectic, huang_mishra_2024, mishra2026majorization, kamat2024simultaneous}, operator theory \cite{bhat2019real, john2022interlacing, kumar2024approximating}, and quantum physics \cite{adesso2004extremal, chen2005gaussian, idel, nicacio2021williamson, hsiang2022entanglement}.

\subsection{Symplectic space and Williamson's theorem}
    A skew-symmetric and non-degenerate bilinear form on a real vector space is called a symplectic form on the vector space. 
    A real vector space with a symplectic form on it is called a symplectic space \footnote{Hermann Weyl \cite{weyl} introduced the term \textit{symplectic} calqued on Greek \textit{sym-plektikos} to mean something similar to \textit{complex}. \textit{Complex} comes from the Latin \textit{com-plexus}, meaning \textit{braided together} (co- + plexus), while \emph{symplectic} comes from the corresponding Greek \textit{sym-plektikos} $(\sigma\upsilon\mu\pi\lambda\epsilon\kappa\tau\iota\kappa$\'{o}$\zeta)$. In both the cases, the part of a word responsible for its lexical meaning comes from the Indo-European root $^*$\textit{ple}$\invbreve{\textit{k}}$-.}, and it is denoted by the pair $\left(\mathscr{V}, \omega \right)$.
    It is well-known that a symplectic space is even dimensional \cite[Proposition~21.1.2]{HormL}.
    Suppose $\mathscr{V}$ is a $2n$-dimensional symplectic space with a symplectic form $\omega$ on it.
    A linear operator $M: \mathscr{V} \to \mathscr{V}$ is said to be symplectic if it preserves the symplectic form, i.e., $\omega(Mu, Mv)=\omega(u, v)$ for all $u, v \in \mathscr{V}$.
    A basis $\{p_1,\ldots, p_n, q_1, \ldots, q_n\}$ of $\mathscr{V}$ is called a \emph{symplectic basis} if it satisfies for all $i,j \in \{1,\ldots, n\}$,
    \begin{align}
        \omega(p_i, p_j)=\omega(q_i, q_j)=0, \quad \omega(p_i, q_j)=\delta_{ij},
    \end{align}
    where $(i,j) \mapsto \delta_{ij}$ is the Kronecker delta function.
    A fundamental result in symplectic linear algebra, known as Williamson's theorem \cite{williamson1936algebraic}, states that if $Q$ is a positive definite quadratic form on $\mathscr{V}$ then there exists a symplectic basis $\{p_1,\ldots, p_n, q_1, \ldots, q_n\}$ of $\mathscr{V}$, and positive numbers $\mu_1,\ldots, \mu_n$ such that for all $(x_1,\ldots, x_n, y_1,\ldots, y_n) \in \mathds{R}^{2n}$,
    \begin{align} \label{eq:general_symp_normal_form}
        Q\left(\sum_{i=1}^n (x_i p_i + y_i q_i)\right)
            = \sum_{i=1}^n \mu_i \left(x_i^2 + y_i^2 \right).
    \end{align}
    We call the diagonalization \eqref{eq:general_symp_normal_form} \emph{Williamson's normal form} of $Q$.

    Our paper is written, without loss of generality, in the language of matrices suitable for the \emph{standard symplectic space}
    $\mathds{R}^{2n}$ equipped with the symplectic form:
    \begin{align}
        \mathds{R}^{2n} \times \mathds{R}^{2n} \ni (x, y) \mapsto x^\top J_{2n} y,
    \end{align}
    where 
        $J_{2n} \coloneqq
            \begin{psmallmatrix}
             0 & I_n \\
             -I_n & 0
            \end{psmallmatrix}
        $,
    $I_n$ being the identity matrix of size $n$.
     We shall drop the subscript $2n$ from $J_{2n}$, and use the notation $J$ instead, when the size of the matrix is clear from the context.
    We will provide interpretations of some of the results for quadratic forms over general symplectic spaces in Section~\ref{sec:interpretation}.

    Symplectic maps on the standard symplectic space are given by \emph{symplectic matrices}, which are $2n \times 2n$ real matrices $M$ that satisfy $M^{\top}JM = J$.
    Positive definite quadratic forms on $\mathds{R}^{2n}$ correspond to $2n \times 2n$ real symmetric positive definite matrices.
    Williamson's theorem states that for every $2n \times 2n$ real symmetric positive definite matrix $A$, there exists a symplectic matrix $M$ such that
    \begin{align}\label{eq:Williamson_decomposition}
        M^{\top} A M 
            &=
            \begin{pmatrix}
             D & 0 \\
             0 & D
            \end{pmatrix},
    \end{align}
    where $D$ is an $n \times n$ diagonal matrix with unique positive diagonal entries (up to ordering), called the symplectic eigenvalues of $A$.
    Several elementary proofs of Williamson's theorem are available in the literature.
    See \cite{folland1989harmonic, simon1999congruences, ikramov2018symplectic}.
    
\subsection{Literature review}
    In his original work \cite{williamson1936algebraic}, Williamson showed that for any $2n \times 2n$ real symmetric matrix $A$ there exists a symplectic matrix $M$ such that $M^{\top}AM$ is a (non-diagonal) sparse matrix. 
    In general, $M^{\top} AM$ may not be a diagonal matrix for any symplectic matrix $M$ much less a diagonal matrix of the form $D \oplus D$ for some $n \times n$ diagonal matrix $D$.
    See the corollary of Theorem~2 in \cite{williamson1936algebraic}.
    Interestingly, if $A$ is positive definite, then it is congruent to a diagonal matrix via a symplectic matrix as stated in \eqref{eq:Williamson_decomposition}.

    Williamson's theorem is known to be generalized to $2n \times 2n$ real symmetric positive semidefinite matrices whose kernels are \emph{symplectic subspaces} of $\mathds{R}^{2n}$.
    More specifically, for a $2n \times 2n$ real symmetric positive semidefinite matrix $A$ there exists a symplectic matrix $M$ such that $M^{\top} AM=D \oplus D$ for some $n \times n$ diagonal matrix $D$ with non-negative diagonal entries if and only if the kernel of $A$ is a symplectic subspace of $\mathds{R}^{2n}$.
    This was stated in \cite[Remark~2.6]{jm}, and explicitly proved in \cite[Theorem~1.3.5]{mishra2021differential}.
    Also, a constructive proof of this extension was recently given in \cite{son2022symplectic}.
    Cruz and Faßbender \cite{cruz_fassbender_2016} established simple algebraic conditions on $2n \times 2n$ complex matrices that are diagonalizable by symplectic equivalence, similarity, or congruence.
    In particular, Theorem~$21$ of \cite{cruz_fassbender_2016} states that for a $2n \times 2n$ (complex) matrix $A$ there exists a (complex) symplectic matrix $M$ such that $M^{\top}AM$ is a diagonal matrix if and only if $A$ is symmetric and $AJ^{\top} A J$ is diagonalizable.
    
    To the best of our knowledge, no precise condition is known for $2n \times 2n$ real symmetric matrices to be diagonalizable in the sense of Williamson's theorem.
    The main aim of this work is to fill this gap.

\subsection{Main contributions}
    In this paper, we establish explicit necessary and sufficient conditions on $2n \times 2n$ real symmetric matrices to be diagonalizable in the sense of Williamson's theorem, and also investigate several implications of it.
    \begin{itemize}
        \item We show that for a $2n \times 2n$ real symmetric matrix $A$ there exists a symplectic matrix $M$ such that $M^{\top} AM=D \oplus D$ where $D$ is an $n \times n$ real diagonal matrix (unique up to ordering of its diagonal entries) if and only if there exist symplectic subspaces $\mathscr{W}_-$, $\mathscr{W}_0$, $\mathscr{W}_+$ of $\mathds{R}^{2n}$ with dimensions $\nu(A), \xi(A), \pi(A)$, respectively such that
        \begin{itemize}
            \item [$\circ$] $\mathscr{W}_-$, $\mathscr{W}_0$, $\mathscr{W}_+$ are pairwise \emph{symplectically orthogonal} to each other
            \item [$\circ$] these subspaces are invariant under $JA$,
            \item [$\circ$] $A$ is negative definite on $\mathscr{W}_-$, the kernel of $A$ is $\mathscr{W}_0$, and $A$ is positive definite on $\mathscr{W}_+$.
        \end{itemize}
        Here $\nu(A), \xi(A), \pi(A)$ denote the number of negative eigenvalues, zero eigenvalues, positive eigenvalues, respectively.
        See Theorem~\ref{thm:williamson_general}.
    \item We introduce a symplectic analog of orthogonal projection, called \emph{symplectic orthogonal projection}, in Definition~\ref{def:sym-orth-proj}, and discuss some properties of it.
    Symplectic orthogonal projections can be of independent interest in symplectic geometry.
    We then re-state the aforementioned result, Theorem~\ref{thm:williamson_general}, in terms of symplectic orthogonal projection.
    See Proposition~\ref{prop:williamson-alternate-sym-projection}.
    This then leads to a more explicit description of the diagonal form in the generalized Williamson's theorem. See Proposition~\ref{prop:sym-eigenvalues-eigenvalues-relation}.
    \item We construct explicit Williamson's decomposition and establish perturbation bounds for the diagonal form for a class of $2n \times 2n$ real symmetric matrices.
    This class, denote by $\operatorname{EigSpSm}(2n)$, consists of $2n \times 2n$ real symmetric matrices whose eigenspaces corresponding to negative eigenvalues, zero eigenvalues, and positive eigenvalues form symplectic subspaces of $\mathds{R}^{2n}$ satisfying the three conditions mentioned above. 
    In particular, $\operatorname{EigSpSm}(2n)$ contains the set of $2n \times 2n$ real positive semidefinite matrices with symplectic kernel.
    The perturbation bounds we obtain generalize known perturbation bounds on symplectic eigenvalues of positive definite matrices given by Bhatia and Jain \cite[Theorem~6]{bhatia2015symplectic}.
    See Section~\ref{sec:special-class}.

    \item We also provide interpretations of the symplectic orthogonal projection and some of the results for quadratic forms in general symplectic spaces in Section~\ref{sec:interpretation} in a coordinate-free fashion, highlighting their geometrical meanings.
    
    \end{itemize}

\subsection{Paper organization}
    We review some basic theory of matrices, linear algebra, and symplectic linear algebra in Section~\ref{sec:review-misc}: Section~\ref{sec:matrices} contains useful concepts from matrix analysis; Section~\ref{sec:lin_algebra_rn} recalls basic theory of subspaces of the Euclidean space $\mathds{R}^n$;
    Section~\ref{sec:review_sp_operations} revisits some basic theory of standard symplectic space $\mathds{R}^{2n}$, and establishes some symplectic operations that are useful for the development of the paper.
    
    We state and prove the main result in Section~\ref{sec:williamson-generalization} (Theorem~\ref{thm:williamson_general}) along with an interesting corollary (Corollary~\ref{cor:williamson_semidefinite}).
    In Section~\ref{sec:williamson-via-sym-proj}, we introduce a symplectic analog of the well-known orthogonal projection called symplectic orthogonal projection (Definition~\ref{def:sym-orth-proj}), and re-state the main result in terms of the symplectic orthogonal projection (Proposition~\ref{prop:williamson-alternate-sym-projection}).
    
     We study Williamson's normal form for a subset of symmetric matrices  $\operatorname{EigSpSm}(2n)$ in Section~\ref{sec:special-class}.
     Here, we explicitly describe the symplectic eigenvalues of matrices in $\operatorname{EigSpSm}(2n)$ (Section~\ref{sec:williamson-construction-1}), construct symplectic matrices achieving the Williamson's normal form (Section~\ref{sec:williamson-construction-2}), and provide perturbation bounds on the symplectic eigenvalues of these matrices (Section~\ref{sec:perturbation_bounds}). 
     Lastly, we provide interpretations of the symplectic orthogonal projection and some of the results for quadratic forms on general symplectic spaces in Section~\ref{sec:interpretation}.

%-----------------------------%-----------------------------%-----------------------------%-----------------------------
    \begin{center}
    \begin{table}[H]
    \caption{Summary of notations and their mathematical definitions.}
    \label{tab:notation}
    \footnotesize
    \begin{tabular}
    [c]{c|l|l}
    \toprule
    Symbol & Meaning & Definition \\\midrule
    $\operatorname{M}(n, k)$ & set of $n \times k$ real matrices &  \\
    $\operatorname{M}(n)$ & set of $n \times n$ real matrices & $\operatorname{M}(n, n)$ \\
    $\operatorname{S}(n)$ & set of symmetric matrices & $\{A \in \operatorname{M}(n)\colon A^{\top}=A\}$ \\
    $\operatorname{Psd}(n)$ & set of positive semidefinite matrices & $\{A \in \operatorname{S}(n)\colon x^{\top} \!A x \geq 0 \forall x \in \mathds{R}^n \}$ \\
    $\operatorname{Pd}(n)$ & set of positive definite matrices  & $\{A \in \operatorname{S}(n)\colon x^{\top} \! A x > 0 \forall x \in \mathds{R}^n\backslash{\{0\}} \}$ \\
    $I_n$ or $I$ & identity matrix of size $n$ \\
    $\operatorname{O}(n)$ & orthogonal group & $\{U \in \operatorname{M}(n)\colon U^{\top}U=I_n\}$ \\
    \bottomrule
    $J_{2n}$ or $J$ & standard symplectic matrix & $\begin{psmallmatrix}
        0 & I_n \\ -I_n & 0
    \end{psmallmatrix}$ \\
    $\operatorname{Sp}(2n, 2k)$ & & $\{M \in \operatorname{M}(2n, 2k)\colon M^{\top}J_{2n} M=J_{2k} \}$ \\
    $\operatorname{Sp}(2n)$ & real symplectic group & $\operatorname{Sp}(2n, 2n)$ \\
    $\operatorname{OSp}(2n)$ & real orthosymplectic group & $\operatorname{O}(2n) \cap \operatorname{Sp}(2n)$
    \\
    $\operatorname{SpS}(2n)$ &  & defined after Remark~\ref{rem:computing-sym-eigenvalues} \\
    $\operatorname{SpPsd}(2n)$ & & $\{A \in \operatorname{Psd}(2n): \ker(A) \cap \ker(A)^{\perp_{\operatorname{s}}}= \{0\}\}$\\
    $\operatorname{EigSpS}(2n)$ & & defined in Section~\ref{sec:special-class}\\
    \bottomrule
    \end{tabular}
    \end{table}
    \end{center}

\section{Review and miscellanea}
\label{sec:review-misc}
    In this section, we establish some notations, and briefly recall some basic concepts from matrix analysis, linear algebra, and symplectic linear algebra.
    We refer the reader to \cite{ma_bhatia, horn2012matrix} for a comprehensive account of theory of matrices, \cite{Johnston_LA_MA} for linear algebra, and \cite{folland1989harmonic, degosson} for symplectic linear algebra.
    A summary of notations with mathematical definitions is provided in Table~\ref{tab:notation}.

\subsection{Matrices}
\label{sec:matrices}
    Let $\operatorname{M}(n, k)$ denote the set of $n \times k$ real matrices.
    We use the shorthand $\operatorname{M}(n)$ for $\operatorname{M}(n,n)$.
    We denote by $\operatorname{S}(n)$ the subset of $\operatorname{M}(n)$ consisting of symmetric matrices.
    For $A \in \operatorname{S}(n)$, we shall use the notations $\nu(A), \xi(A)$, $\pi(A)$ to denote the number of negative eigenvalues, zero eigenvalues, positive eigenvalues of $A$, respectively.
    If $K \in \operatorname{M}(n)$ is an invertible matrix then the Sylvester's law of inertia states that for any $A \in \operatorname{Sm}(n)$, we have $\nu(A)=\nu(K^{\top}AK)$, $\xi(A)=\xi(K^{\top}AK)$, and $\pi(A)=\pi(K^{\top}AK)$. 
    See \cite[Theorem~4.5.8]{horn2012matrix}.

    We denote by $\operatorname{Psd}(n)$ and $\operatorname{Pd}(n)$ the subsets of $\operatorname{S}(n)$ consisting of positive semidefinite and positive definite matrices, respectively. 
    Let $\operatorname{O}(n)$ denote the real orthogonal group in dimension $n$.
    A matrix $A \in \operatorname{M}(n)$ is called normal if $A^{\top}A=AA^{\top}$.
    For every $B \in \operatorname{Psd}(n)$, there exists a unique $B^{1/2} \in \operatorname{Psd}(n)$ such that
    $(B^{1/2})^2 = B$.
    The matrices $B$ and $B^{1/2}$ have the same range, and hence the same rank.
    See \cite[Theorem~7.2.6]{horn2012matrix}.
    Every symmetric matrix $C \in \operatorname{S}(n)$ can be expressed as a difference of two positive semidefinite matrices $C=C_+-C_-$, where
    \begin{align}
        C_- &\coloneqq \dfrac{1}{2}(|C|-C), \label{eq:positive_part_matrix}\\
        C_+ &\coloneqq \dfrac{1}{2}(|C|+C), \label{eq:negative_part_matrix}
    \end{align}
    and $|C| \coloneqq (C^2)^{1/2}$ is the absolute value of $C$.
    We have $\operatorname{rank}(C)=\operatorname{rank}(C_+)+\operatorname{rank}(C_-)$ and $C_+ C_- = C_- C_+=0$.
    See Proposition~$4.1.13$ of \cite{horn2012matrix}.

\subsection{Linear algebra on $\mathds{R}^{n}$}\label{sec:lin_algebra_rn}
    We denote by $\langle \cdot , \cdot \rangle$ the Euclidean inner product given for all $x, y \in \mathds{R}^n$ by $\langle x, y \rangle \coloneqq x^{\top}y$.
    Let $\mathscr{W}$ be a linear subspace of $\mathds{R}^{n}$.
    $\mathscr{W}$ is said to be an invariant subspace of $A \in \operatorname{M}(n)$ if for all $w \in \mathscr{W}$,  $Aw \in \mathscr{W}$.
    We say that $A$ is positive definite on $\mathscr{W}$ if $\langle w, Aw \rangle > 0$ for all non-zero $w \in \mathscr{W}$.
    We say $A$ is negative definite on $\mathscr{W}$ if $-A$ is positive-definite on $\mathscr{W}$.
    The orthogonal complement of $\mathscr{W}$ is defined as
    \begin{align}
        \mathscr{W}^{\perp}\coloneqq \{u \in \mathds{R}^n: \langle u, w \rangle =0, \ \forall w \in \mathscr{W}\}.
    \end{align}
    %We have $\mathds{R}^n = \mathscr{W}\oplus \mathscr{W}^{\perp}$, which means that for every $x \in \mathds{R}^{n}$, there exist unique $w \in \mathscr{W}$ and $w' \in \mathscr{W}^{\perp}$ such that $x=w+w'$.
    A matrix $P \in \operatorname{S}(n)$ is called an orthogonal projection onto $\mathscr{W}$ if $Pw=w$ and $Pw'=0$ for all $w \in \mathscr{W}$ and $w' \in \mathscr{W}^{\perp}$.
    Any matrix $Q \in \operatorname{S}(n)$ that satisfies $Q^2=Q$ is an orthogonal projection onto $\operatorname{range}(Q)$.

    % We will use the following version of the well-known spectral theorem of symmetric matrices.
    % A standard proof of this result uses block diagonal representation of a symmetric matrix in an orthonormal basis extending any orthonormal basis of its invariant subspace.
    % See Theorem~6.1.6 of \cite{watkins2004fundamentals}.
    % \begin{boxed}
    %     \begin{lemma}\label{lem:invariant_subspace_eigenspace}
    %     Let $A \in \operatorname{S}(n)$ and $\mathscr{W} \subseteq \mathds{R}^n$ be an invariant subspace of $A$. 
    %     Then $\mathscr{W}$ is spanned by some eigenvectors of $A$.
    %     If the restriction of $A$ on $\mathscr{W}$ is positive definite, then $\mathscr{W}$ is spanned by some eigenvectors of $A$ corresponding to its positive eigenvalues.        
    % \end{lemma}
    % \end{boxed}

%-----------------------------%-----------------------------%-----------------------------%-----------------------------

\subsection{Symplectic linear algebra on $\mathds{R}^{2n}$}\label{sec:review_sp_operations}

    The {\it symplectic orthogonal complement} of a subset $\mathscr{X} \subseteq \mathds{R}^{2n}$ is defined as
    \ba
    {
    	\mathscr{X}^{\perp_{\operatorname{s}}} \coloneqq \{u \in \mathds{R}^{2n}: \forall v \in \mathscr{X}, \langle u, Jv \rangle = 0 \}.
    }
    
    A linear subspace $\mathscr{W}$ of $\mathds{R}^{2n}$ is called a symplectic subspace if for every $u \in \mathscr{W}$ there exists $v \in \mathscr{W}$ such that $\langle u, Jv \rangle \neq 0$.
    By definition, $\mathscr{W}$ is a symplectic subspace of $\mathds{R}^{2n}$ if and only if
    $\mathscr{W} \cap \mathscr{W}^{\perp_{\operatorname{s}}}=\{0\}$.
    Let $\mathscr{W}$ be a symplectic subspace of $\mathds{R}^{2n}$.
    Then $\mathscr{W}$ has even dimension, say $2k$, and it has a {\it symplectic basis} $\{u_1,\ldots, u_k, v_1,\ldots, v_k\}$ that satisfies for all $1\leq i,j \leq n$:
    \begin{align}
        \langle u_i, J v_j \rangle &= \delta_{ij}, \\
        \langle u_i, J u_j \rangle &= 0, \\
        \langle v_i, J v_j \rangle &= 0.
    \end{align}
    Here $\delta_{ij}=0$ if $i\neq j$ and $\delta_{ij}=1$ if $i=j$. 
    We have $\dim(\mathscr{W})+\dim(\mathscr{W}^{\perp_{\operatorname{s}}})=2n$ and $\Par{\mathscr{W}^{\perp_{\operatorname{s}}}}^{\perp_{\operatorname{s}}}= \mathscr{W}$.
    See \cite[Section~1.2]{degosson}.
    We say that two symplectic subspaces $\mathscr{W}$ and $\mathscr{Z}$ are said to be symplectically orthogonal to each other if $\mathscr{Z} \subseteq \mathscr{W}^{\perp_{\operatorname{s}}}$.
    
    Let $\operatorname{Sp}(2n, 2k)$ denote the set of $2n \times 2k$ real  matrices $M$ that satisfy $M^{\top} J_{2n}M=J_{2k}$.
    We use the shorthand $\operatorname{Sp}(2n)$ for $\operatorname{Sp}(2n, 2n)$.
    The set $\operatorname{Sp}(2n)$ consists of $2n \times 2n$ real symplectic matrices, and it is known as the symplectic group.
    For every $M \in \operatorname{Sp}(2n, 2k)$, $\operatorname{range}(M)$ is a symplectic subspace of $\mathds{R}^{2n}$, and the columns of $M$ form a symplectic basis of $\operatorname{range}(M)$.
    See \cite[Section~1.2.1]{degosson}.
    We denote by $\OrSp(2n) \coloneqq \Or(2n) \cap \Sp(2n)$ the set of {\it orthosymplectic} matrices.

    Let $n_1,\ldots, n_k$ be positive integers, and $X_i \in \operatorname{M}(n_i)$ for $1\leq i \leq k$.
    Denote by $\oplus X_i$ the usual direct sum of the matrices $X_1,\ldots, X_k$.
    Suppose $A_i \in \operatorname{M}(2n_i)$ is partitioned into blocks as 
    \begin{align}
        A_i = 
            \begin{pmatrix}
                E_i & F_i \\
                G_i & H_i
            \end{pmatrix},
    \end{align}
    where $E_i, F_i, G_i, H_i \in \operatorname{M}(n_i)$ for all $1\leq i \leq k$.
    The {\it $\operatorname{s}$-direct sum} of $A_1,\ldots, A_k$ is defined by
    \begin{align}
        \oplus^{\operatorname{s}} A_i 
            \coloneqq
                \begin{pmatrix}
                \oplus E_i & \oplus F_i \\
                \oplus G_i & \oplus H_i
            \end{pmatrix}.            
    \end{align}

    Let $M$ and $N$ be $2n \times 2k$ and $2n \times 2\ell$ matrices whose columns are $u_1,\ldots, u_{k}, v_1,\ldots,v_k$ and $x_1,\ldots, x_{\ell}, y_1,\ldots,y_\ell$, respectively. 
    Define the {\it symplectic concatenation} of $M$ and $N$ to be the following $2n \times 2(k+\ell)$ matrix given by
    \begin{align}
        M \diamond N \coloneqq \Brac{u_1,\ldots, u_{k},x_1,\ldots, x_{\ell},v_1,\ldots,v_k, y_1,\ldots,y_\ell }.
    \end{align}

%-----------------------------%-----------------------------%-----------------------------%-----------------------------
\section{Williamson's theorem for symmetric matrices}\label{sec:williamson-generalization}
    Generalizing Williamson's theorem to symmetric matrices is the main objective of this section.
    We begin by building some intuition towards generalization of the theorem.
    Let $A \in \operatorname{S}(2n)$ for which
    there exists $M \in \operatorname{Sp}(2n)$ such that
    \begin{align}\label{eq:extended_williamson_alternate_intuition}
        M^{\top} A M 
            &=
            \begin{pmatrix}
             D & 0 \\
             0 & D
            \end{pmatrix},
    \end{align}
    where $D$ is an $n \times n$ diagonal matrix.
    We shall refer to \eqref{eq:extended_williamson_alternate_intuition} as a Williamson's decomposition of $A$.
    Since the symplectic matrix $M$ satisfies $M^{-T}=J M J^{\top}$, \eqref{eq:extended_williamson_alternate_intuition} gives
    \begin{align}\label{eq:extended_williamson_alternate_intermediate}
        AM= JMJ^{\top}\begin{pmatrix}
             D & 0 \\
             0 & D
            \end{pmatrix}.
    \end{align}
    Let $u_1,\ldots, u_n, v_1,\ldots, v_n$ denote the columns of $M$ and $d_1,\ldots, d_n$ denote the diagonal elements of $D$.
    Then \eqref{eq:extended_williamson_alternate_intermediate} implies for all $1\leq i \leq n$:
    \begin{align}
        A u_i &= d_i J v_i, \label{eq:extended_williamson_alternate_1}\\
        A v_i &= -d_i J u_i. \label{eq:extended_williamson_alternate_2}
    \end{align} 
    Define index sets:
    \begin{align}
        \mathcal{I}_- &\coloneqq \{i: 1 \leq i \leq n, d_i<0\},\\
        \mathcal{I}_0 &\coloneqq \{i: 1 \leq i \leq n, d_i=0\},\\
        \mathcal{I}_+ &\coloneqq \{i: 1 \leq i \leq n, d_i>0\},
    \end{align}
    and subspaces:
    \begin{align}
        \mathscr{W}_- &\coloneqq \operatorname{span}\{u_i, v_i: i \in \mathcal{I}_-\},\\
        \mathscr{W}_0 &\coloneqq \operatorname{span}\{u_i, v_i: i \in \mathcal{I}_0\},\\
        \mathscr{W}_+ &\coloneqq \operatorname{span}\{u_i, v_i: i \in \mathcal{I}_+\}.
    \end{align}
    By construction, $\mathscr{W}_-, \mathscr{W}_0, \mathscr{W}_+$ are symplectic subspaces and are pairwise symplectically orthogonal to each other.
    Also, by the Sylvester's law of inertia, we have $\dim(\mathscr{W}_-)=\nu(A)$, $\dim(\mathscr{W}_0)=\xi(A)$, and $\dim(\mathscr{W}_+)=\pi(A)$ so that the dimensions of these subspaces add to $2n$.
    The relations \eqref{eq:extended_williamson_alternate_1} and \eqref{eq:extended_williamson_alternate_2} imply that these subspaces are invariant under $JA$.
    It is also easy to verify that $A$ is negative definite on $\mathscr{W}_-$.
    Indeed, let $x \in \mathscr{W}_-$ be any non-zero vector given by $x=\sum_{i \in \mathcal{I}_-} (a_i u_{i}+b_i v_{i})$, where $a_i, b_i \in \mathds{R}$ for all $i \in \mathcal{I}_-$.
    We have
    \begin{align}
        \langle x, Ax \rangle 
            &= \bigg\langle \sum_{i \in \mathcal{I}_-} (a_i u_{i}+b_i v_{i}), \sum_{j \in \mathcal{I}_-} (a_{j}A u_{j}+b_{j} Av_{j}) \bigg \rangle \\
            &= \sum_{i, j \in \mathcal{I}_-}   \langle a_i u_{i}+b_i v_{i},  a_{j}d_{j}J v_{j}-b_{j} d_{j}Ju_{j}  \rangle \\
            &= \sum_{i, j \in \mathcal{I}_-} 
            \left(  a_{i}a_{j} d_{j} \langle  u_{i}, J v_{j} \rangle -  a_ib_{j} d_{j} \langle u_i, Ju_{j} \rangle \right. \notag \\
            &\hspace{2cm} \left.+  b_{i}a_{j} d_{j} \langle  v_{i}, J v_{j} \rangle -  b_ib_{j} d_{j} \langle v_i, Ju_{j} \rangle \right) \\
            &= \sum_{i \in \mathcal{I}_-} d_i (a_i^2 + b_i^2) \\
            &<0.
    \end{align}
    The last inequality follows from the fact that $d_i < 0$ for all $i \in \mathcal{I}_-$.
    A similar argument shows that $A$ is positive definite on $\mathscr{W}_+$.
    Also, we obviously have $\operatorname{ker}(A)=\mathscr{W}_0$.
    
    To summarise everything, the following are necessary conditions on any $A \in \operatorname{S}(2n)$ that is diagonalizable in the sense of Williamson's theorem:
    \begin{itemize}
        \item [] $\mathbf{Condition~(i)}$  There exist pairwise symplectically orthogonal symplectic subspaces $\mathscr{W}_-, \mathscr{W}_0, \mathscr{W}_+$ with dimensions $\nu(A), \xi(A), \pi(A)$, respectively.
        \item [] $\mathbf{Condition~(ii)}$  Each of these symplectic subspaces is invariant under $JA$.
        \item [] $\mathbf{Condition~(iii)}$ $A$ is negative definite on $\mathscr{W}_-$, the kernel of $A$ is $\mathscr{W}_0$,  and $A$ is positive definite on $\mathscr{W}_+$.
    \end{itemize}

    In the following theorem, we show that the above three conditions are sufficient for a symmetric matrix to be diagonalizable in the sense of Williamson's theorem.
    \begin{boxed}
    \begin{theorem}\label{thm:williamson_general}
    For $A \in \operatorname{S}(2n)$ there exists $M \in \Sp(2n)$ and an $n \times n$ diagonal matrix $D$ such that
    \begin{align}\label{eq:williamson-general-statement}
        M^{\top} A M =  \begin{pmatrix}
             D & 0 \\
             0 & D
            \end{pmatrix}
    \end{align}
    if and only $A$ satisfies $\mathbf{Condition~(i)}, \mathbf{Condition~(ii)}$, and $\mathbf{Condition~(iii)}$.
    The diagonal matrix $D$ so obtained is unique up to ordering of its diagonal entries.
    Moreover, the diagonal entries of $D$ and $-D$ combined together constitute the $2n$ eigenvalues of $i  J A$, where $i=\sqrt{-1}$.
    \end{theorem}
    \end{boxed}
    \begin{proof}
        The necessity of the given conditions is already established in the beginning of the section.
        In what follows, we give an argument for sufficiency of these conditions.
        
        Let $A \in \operatorname{S}(2n)$ and suppose $\mathscr{W}_-, \mathscr{W}_0, \mathscr{W}_+$ are symplectic subspaces of $\mathds{R}^{2n}$ that satisfy $\mathbf{Condition~(i)}, \mathbf{Condition~(ii)}$, and $\mathbf{Condition~(iii)}$ for  $A$. 
        Let $k=\frac{1}{2}\nu(A), \ell=\frac{1}{2}\xi(A)$, and $m=\frac{1}{2}\pi(A)$.
        Choose $M_- \in \operatorname{Sp}(2n, 2k)$, $M_0 \in \operatorname{Sp}(2n, 2\ell)$, and $M_+ \in \operatorname{Sp}(2n, 2m)$
        such that $\operatorname{range}(M_-) = \mathscr{W}_-$, $\operatorname{range}(M_0) = \mathscr{W}_0$, and $\operatorname{range}(M_+) = \mathscr{W}_+$.
        By $\mathbf{Condition~(i)}$ we have $M_- \diamond M_0 \diamond M_+ \in \operatorname{Sp}(2n)$.
        $\mathbf{Condition~(iii)}$ implies $-M_-^{\top} A M_- \in \operatorname{Pd}(2k)$ and $M_+^{\top} A M_+ \in \operatorname{Pd}(2m)$. 
        By Williamson's theorem, we thus get $Q_- \in \operatorname{Sp}(2k)$ and $Q_+ \in \operatorname{Sp}(2m)$ such that
        \begin{align}
            Q_-^{\top}M_-^{\top} A M_- Q_- 
                &= 
                \begin{pmatrix}
                    D_- & 0 \\
                    0 & D_-
                \end{pmatrix}, \\
            Q_+^{\top}M_+^{\top} A M_+ Q_+ 
                &= 
                \begin{pmatrix}
                    D_+ & 0 \\
                    0 & D_+
                \end{pmatrix},
        \end{align}
        where $D_- < 0$ and $D_+ > 0$ are diagonal matrices of size $k \times k$ and $m \times m$, respectively.
        Set $M \coloneqq \left(M_- Q_- \right) \diamond M_0 \diamond \left(M_+ Q_+ \right)$.
        It is easy to check that $M \in \operatorname{Sp}(2n)$.
        In what follows, we show that $M$ diagonalizes $A$ in the sense of Williamson's theorem.

        By $\mathbf{Condition~(ii)}$, the columns of $JAM_+$ lie in the subspace $\mathscr{W}_+$.
        Since $\mathscr{W}_-$ and $\mathscr{W}_+$ are symplectically orthogonal to each other, we have $M_-^{\top} J JAM_+=0$  implying $M_-^{\top} AM_+=0$.
        Also, we have $AM_0=0$ which implies that $M_-^{\top}AM_0=0$ and $M_+^{\top}AM_0=0$.
        Therefore, we get
        \begin{align}
            M^{\top} A M 
                &= \left[\left(M_- Q_- \right) \diamond M_0 \diamond \left(M_+ Q_+ \right)\right]^{\top} A \left[\left(M_- Q_- \right) \diamond M_0 \diamond \left(M_+ Q_+ \right) \right] \\
                &= \left[\left(M_- Q_- \right) \diamond M_0 \diamond \left(M_+ Q_+ \right)\right]^{\top}  \left[\left(A M_- Q_- \right) \diamond \Par{A M_0} \diamond \left(AM_+ Q_+ \right) \right] \\
                &= \left(Q_-^{\top} M_-^{\top} A M_- Q_- \right) \oplus_{\operatorname{s}} \left( M_0^{\top} A M_0 \right)
                \oplus_{\operatorname{s}} \left(Q_+^{\top} M_+^{\top} A M_+ Q_+ \right) \\
                &= \begin{pmatrix}
                    D_- & 0 \\
                    0 & D_-
                \end{pmatrix} \oplus_{\operatorname{s}}
                \begin{pmatrix}
                    0_\ell & 0_\ell \\
                    0_\ell & 0_\ell
                \end{pmatrix}
                \oplus_{\operatorname{s}}
                \begin{pmatrix}
                    D_+ & 0 \\
                    0 & D_+
                \end{pmatrix} \\
                &= \begin{pmatrix}
                    D & 0 \\
                    0 & D
                \end{pmatrix},
        \end{align}
        where $D\coloneqq D_- \oplus 0_\ell \oplus D_+$ and $0_\ell$ denotes the zero matrix of size $\ell \times \ell$.

        The uniqueness of the diagonal form $D$ and the fact that the combined diagonal entries of $D$ and $-D$ form the eigenvalues of $i  JA$ are established by Pereira et al.~\cite[Section~5]{pereira2021symplectic}.
    \end{proof}
    \begin{boxed}
        \begin{remark}\label{rem:computing-sym-eigenvalues}
    	\emph{Pereira et al.~\cite{pereira2021symplectic} provided a method of computing Williamson's decomposition of a $2n \times 2n$ (complex) symmetric matrix $A$, given that $A$ is guaranteed to admit such a decomposition.
        Theorem~\ref{thm:williamson_general} complements the work of Pereira et al.~\cite{pereira2021symplectic} in the sense that it provides a characterization of $A$ for existence of its Williamson's decomposition.}
        \end{remark}
    \end{boxed}

    Let  $\SpSm(2n)$ denote the subset of $\operatorname{Sm}(2n)$ consisting of matrices satisfying $\mathbf{Condition~(i)}$, $\mathbf{Condition~(ii)}$, and $\mathbf{Condition~(iii)}$.
    In view of Theorem~\ref{thm:williamson_general}, for every $A \in \SpSm(2n)$, there exists $M \in \Sp(2n)$ and a unique $n \times n$ diagonal matrix $D$ with diagonal diagonal entries in ascending order such that $M^{\top} A M = D \oplus D$.
    We refer to the diagonal elements of $D$ as the symplectic eigenvalues of $A$.
    Thus, a matrix in $\SpSm(2n)$ can have negative, zero, or positive symplectic eigenvalues.

    Let $\operatorname{SpPsd}(2n)$ denote the set of $2n \times 2n$ real symmetric positive semidefinite matrices with symplectic kernel.
    As a corollary of Theorem~\ref{thm:williamson_general}, we get the following known result which states that every matrix in $\operatorname{SpPsd}(2n)$ exhibits Williamson's decomposition.
    See \cite[Remark~2.6]{jm} and \cite[Section~2]{son2022symplectic}.

    \begin{boxed}
    \begin{corollary}\label{cor:williamson_semidefinite}
     We have $\operatorname{SpPsd}(2n) \subset \SpSm(2n)$.
    \end{corollary}
    \end{boxed}
    \begin{proof}
        Let $A \in \operatorname{SpPsd}(2n)$.
        Choose $\mathscr{W}_-=0$, $\mathscr{W}_0=\operatorname{ker}(A)$, and $\mathscr{W}_+ = \mathscr{W}_0^{\perp_{\operatorname{s}}}$.
        These symplectic subspaces clearly satisfy $\mathbf{Condition~(i)}$ and $\mathbf{Condition~(iii)}$.
        It is also straightforward to see that $\mathscr{W}_-$ and $\mathscr{W}_0$ are invariant under $JA$.
        It remains to show that $\mathscr{W}_+$ is invariant under $JA$.
       	We have $\mathds{R}^{2n} = \mathscr{W}_0 \oplus \mathscr{W}_+$.
        Let $y \in \mathscr{W}_+$ be arbitrary.
        For any $x \in \mathscr{W}_0$, we have
        \begin{align}
            \langle x, J(J A) y \rangle 
             &= -\langle x,  A y \rangle \\
             &=-\langle A x,   y \rangle \\
             &= -\langle 0,   y \rangle \\
             &=0.
        \end{align}
        This by definition means $JA y \in \mathscr{W}_0^{\perp_{\operatorname{s}}}=\mathscr{W}_+$, implying that $\mathscr{W}_+$ is invariant under $JA$.
        This shows that $\mathbf{Condition~(ii)}$ is also satisfied by $\mathscr{W}_-$, $\mathscr{W}_0$, $\mathscr{W}_+$ for $A$ and hence $A \in \SpSm(2n)$.
    \end{proof}

\section{General Williamson's theorem via symplectic orthogonal projection}
\label{sec:williamson-via-sym-proj}
    In this section we introduce a symplectic analog of orthogonal projection, call it symplectic orthogonal projection, and provide an alternate statement for the general Williamson's theorem in terms of symplectic orthogonal projection.
    
    Let $\mathscr{W}$ be a $2k$ dimensional symplectic subspace of $\mathds{R}^{2n}$.
    Let $M \in \operatorname{Sp}(2n, 2k)$ be any matrix such that $\operatorname{range}(M)=\mathscr{W}$.
    The matrix $P_M \coloneqq JMM^{\top} J^{\top}$ is called the symplectic projection corresponding to $M$.
    It is a positive semidefinite matrix with kernel $\mathscr{W}^{\perp_{\operatorname{s}}}$.
    See \cite[Section~5]{jm}.
    It is known that for $N \in \operatorname{Sp}(2n, 2k)$, the equality $P_M=P_N$ holds if and only there exists $U \in \operatorname{OSp}(2k)$ such that $N=MU$ \cite[Proposition~5.1]{jm}.
    Consequently, we have that $\operatorname{range}(N)=\mathscr{W}$ is a necessary but not a sufficient condition for the symplectic projection $P_N$ to be equal to $P_M$ (for instance, choose $N=MU$ for $U \in \Sp(2k) \backslash \OrSp(2k)$).
    However, it is interesting to observe that the condition $\operatorname{range}(N)=\mathscr{W}$ is necessary and sufficient for the equality $P_N JP_N=P_M JP_M$.
    Moreover, the matrix $J^{\top}P_M JP_M$ restricted to $\mathscr{W}$ is the identity operator and its kernel is $\mathscr{W}^{\perp_{\operatorname{s}}}$ as shown in the following proposition.
    \begin{boxed}
    \begin{proposition}\label{prop:symp_subspace_projection}
        Let $\mathscr{W}$ be a symplectic subspace of $\mathds{R}^{2n}$.
        Let $M \in \operatorname{Sp}(2n, 2k)$ such that $\operatorname{range}(M)=\mathscr{W}$, and let $P_M$ be the symplectic projection corresponding to $M$.
        Then for all $x \in \mathscr{W}$, we have $J^{\top} P_{M}JP_M x = x$. 
        Also, $\operatorname{ker}(J^{\top} P_{M}JP_M)=\mathscr{W}^{\perp_{\operatorname{s}}}$.
    \end{proposition}
    \end{boxed}
    \begin{proof}
        Let $M \in \operatorname{Sp}(2n, 2k)$ such that $\operatorname{range}(M)=\mathscr{W}$, and 
        let $u_1,\ldots, u_k, v_1,\ldots, v_k$ be the columns of $M$.
        Denote by $e_1,\ldots, e_{2n}$ the standard unit vectors of $\mathds{R}^{2n}$.
        For all $1 \leq i \leq r$ we have
        \begin{align}
            P_M u_i 
                &= J M M^{\top} J^{\top} u_i \\
                &= J M e_{i+n} \\
                &= J v_i.
        \end{align}
        Similarly, we get $P_M v_i = -J u_i$.
        These observations give the following: 
        \begin{align}
            J^{\top}P_MJ P_M u_i
                &= J^{\top}P_MJ^2 v_i \\
                &= J P_M v_i \\
                &= -J^2 u_i \\
                &= u_i.
        \end{align}
        A similar argument gives $J^{\top}P_MJ P_M v_i=v_i$.
        Consequently,  for all $x \in \mathscr{W}$, we have $J^{\top}P_MJ P_M x = x$.

        \sloppy We have $\operatorname{ker}(J^{\top}P_MJ P_M ) \supseteq  \operatorname{ker}(P_M) = \mathscr{W}^{\perp_{\operatorname{s}}}$, and $\operatorname{range}(J^{\top}P_MJ P_M) \supseteq \mathscr{W}$.
        The rank-nullity theorem, combined with the fact that $\dim(\mathscr{W})+ \dim(\mathscr{W}^{\perp_{\operatorname{s}}})=2n$, implies that $\operatorname{ker}(J^{\top}P_MJ P_M)=\mathscr{W}^{\perp_{\operatorname{s}}}$.
    \end{proof}
    Proposition~\ref{prop:symp_subspace_projection} states that associated with every symplectic subspace is a \emph{unique} matrix that acts as the identity on the symplectic subspace and its kernel is the symplectic complement of the given symplectic subspace.
    This leads to the following definition of {\it symplectic orthogonal projection} onto a symplectic subspace.
    \begin{boxed}
    \begin{definition}\label{def:sym-orth-proj}
        Let $\mathscr{W}$ be a symplectic subspace of $\mathds{R}^{2n}$.
        We call the $2n \times 2n$ real matrix $\Pi$, given by
        \begin{align}
            \Pi x =
                &\begin{cases}
                    x & \text{if }x \in \mathscr{W}, \\
                    0 & \text{if } x \in \mathscr{W}^{\perp_{\operatorname{s}}},
                \end{cases}
        \end{align}
        the symplectic orthogonal projection onto $\mathscr{W}$.
        It is given by
        \begin{align}
            \Pi = J^{\top} P_M J P_M,
        \end{align}
        for any $M \in \Sp(2n, 2k)$ such that $\operatorname{range}(M)=\mathscr{W}$.
    \end{definition}
    \end{boxed}

    \begin{boxed}
        \begin{remark}
            Symplectic orthogonal projections are precisely $2n \times 2n$ real matrices $\Pi$ that satisfy
                \begin{itemize}
                    \item $\ker(\Pi)$ is a symplectic subspace,
                    \item $\operatorname{range}(\Pi)=\ker(\Pi)^{\perp_{\operatorname{s}}}$,
                    \item $\Pi^2=\Pi$.
                \end{itemize}
        \end{remark}
    \end{boxed}

    \begin{boxed}
    \begin{proposition}\label{prop:transpose-symp-subspace-projection}
        Let $\mathscr{W}$ be a symplectic subspace of $\mathds{R}^{2n}$, and let $\Pi$ be the associated symplectic orthogonal projection.
        Then $\Pi^{\top}$ is the symplectic orthogonal projection onto $J \mathscr{W}$.
    \end{proposition}
    \end{boxed}
    \begin{proof}
        Let $M \in \Sp(2n, 2k)$ such that $\operatorname{range}(M)=\mathscr{W}$.
        We have $P_{JM}= MM^{\top}$.
        This gives
        \begin{align}
            \Pi^{\top} 
                &= \Par{J^{\top} P_M J P_M}^{\top} \\
                &= P_M J^{\top} P_M J \\
                &=  J MM^{\top} J^{\top} J^{\top} J MM^{\top} J^{\top} J \\
                &= J^{\top} P_{JM} J P_{JM}.
        \end{align}
    \end{proof}

    We now state Theorem~\ref{thm:williamson_general} in terms of symplectic orthogonal projections as follows.
    \begin{boxed}
    \begin{proposition}\label{prop:williamson-alternate-sym-projection}
        Let $A \in \operatorname{Sm}(2n)$.
        We have $A \in \operatorname{SpS}(2n)$ if and only if there exist symplectic orthogonal projections $\Pi_-, \Pi_0, \Pi_+$ satisfying the following conditions:
        \begin{enumerate}
            \item [$(i)$] $\Pi_{-}\Pi_{0}=\Pi_{-}\Pi_{+}=\Pi_{0}\Pi_{+}=0$, and $\Pi_{-}+\Pi_{0}+\Pi_{+}=I_{2n}$.
            \item [$(ii)$] $A =  \Pi_{-}^{\top} A \Pi_{-}+\Pi_{+}^{\top} A \Pi_{+} $.
            \item [$(iii)$] $\Pi_{-}^{\top} A \Pi_{-}$ is negative definite on $\operatorname{range}(\Pi_-)$ and $\Pi_{+}^{\top} A \Pi_{+}$ is positive definite on $\operatorname{range}(\Pi_+)$.
        \end{enumerate}
    \end{proposition}
    \end{boxed}
    \begin{proof}
	The ``if'' part is straightforward. 
	Suppose there exist symplectic orthogonal projections $\Pi_-, \Pi_0, \Pi_+$ satisfying the given conditions.
	Choose $\mathscr{W}_{-}=\operatorname{range}(\Pi_-)$, $\mathscr{W}_{0}=\operatorname{range}(\Pi_0)$, and $\mathscr{W}_{+}=\operatorname{range}(\Pi_+)$.
	It is easy to see that the symplectic subspaces  $\mathscr{W}_{-}$, $\mathscr{W}_{0}$, $\mathscr{W}_{+}$ satisfy $\mathbf{Condition~(i)}$, $\mathbf{Condition~(ii)}$, and $\mathbf{Condition~(iii)}$.
	Therefore, we have $A \in \SpSm(2n)$.
	
	We now prove the ``only if'' part.
        Suppose $A \in \SpSm(2n)$.
        Then there exist symplectic subspaces $\mathscr{W}_{-}, \mathscr{W}_0, \mathscr{W}_+$ satisfying $\mathbf{Condition~(i)}$, $\mathbf{Condition~(ii)}$, and $\mathbf{Condition~(iii)}$ for $A$.
        Let $\Pi_-$, $\Pi_0$, and $\Pi_+$ be the symplectic orthogonal projections onto $\mathscr{W}_{-}$, $\mathscr{W}_0$, and $\mathscr{W}_+$, respectively. 
        $\mathbf{Condition~(i)}$ implies that $\Pi_{-}\Pi_{0}=\Pi_{-}\Pi_{+}=\Pi_{0}\Pi_{+}=0$ and $\Pi_{-}+\Pi_{0}+\Pi_{+}=I_{2n}$.
        For any $x \in \mathscr{W}_-$, $y \in \mathscr{W}_0$, and $z \in \mathscr{W}_+$, we get
        \begin{align}
            \Par{\Pi_{-}^{\top} A \Pi_{-}+ \Pi_{+}^{\top} A \Pi_{+}}(x+y+z)
                &=   \Pi_{-}^{\top} A \Pi_{-} x + \Pi_{+}^{\top} A \Pi_{+} z\\
                &= \Pi_{-}^{\top} A x + \Pi_{+}^{\top} A z \\
                &= Ax + Az \label{eq:condition-ii-implication}\\
                &= A(x+y+z).
        \end{align}
        The equality \eqref{eq:condition-ii-implication} follows from Proposition~\ref{prop:transpose-symp-subspace-projection} and the fact that $\mathscr{W}_-$ and $\mathscr{W}_+$ are invariant under $JA$, which is given by $\mathbf{Condition~(ii)}$. 
        We thus have $A=\Pi_{+}^{\top} A \Pi_{+} + \Pi_{-}^{\top} A \Pi_{-}$.
        Lastly, $\Pi_{+}^{\top} A \Pi_{+}$ being positive definite on $\mathscr{W}_+$ and $\Pi_{-}^{\top} A \Pi_{-}$ being negative definite on $\mathscr{W}_-$ follows directly from $\mathbf{Condition~(iii)}$.
    \end{proof}

    We know that the symplectic eigenvalues of $A \in \Pd(2n)$ are the positive eigenvalues of the Hermitian matrix $i  A^{1/2}JA^{1/2}$. 
    We state an analogous fact for matrices in  $\SpSm(2n)$ as follows.        
    \begin{boxed}
        \begin{proposition}\label{prop:sym-eigenvalues-eigenvalues-relation} 
            Let $A \in \SpSm(2n)$, and let $\Pi_-$, $\Pi_0$, and $\Pi_+$ be symplectic orthogonal projections given by Proposition~\ref{prop:williamson-alternate-sym-projection}. 
            Then $A$ has $\frac{1}{2}\rank\Par{\Pi_0}$ zero symplectic eigenvalues.
            The negative symplectic eigenvalues of $A$ are the negative eigenvalues of $i  \Par{-\Pi_-^{\top} A \Pi_-}^{1/2} J \Par{-\Pi_-^{\top} A \Pi_-}^{1/2}$, and the positive symplectic eigenvalues are the positive eigenvalues of $i  \Par{\Pi_+^{\top} A \Pi_+}^{1/2} J \Par{\Pi_+^{\top} A \Pi_+}^{1/2}$.
        \end{proposition}
    \end{boxed}
    \begin{proof}
            We know from Proposition~\ref{prop:williamson-alternate-sym-projection} that $-\Pi_-^{\top} A \Pi_-$ and $\Pi_+^{\top} A \Pi_+$ are positive semidefinite matrices.
            Also, $\Ker(-\Pi_-^{\top} A \Pi_-)=\Ker(\Pi_-)$ and $\Ker(\Pi_+^{\top} A \Pi_+)=\Ker(\Pi_+)$,
            which follow from the facts that $\Pi_-^{\top} A \Pi_-$ is negative definite on $\ran\Par{\Pi_-}$ and $\Pi_+^{\top} A \Pi_+$ is positive definite on $\ran\Par{\Pi_+}$.
            Therefore, Williamson's decompositions of $\Pi_-^{\top} A \Pi_-$ and $\Pi_+^{\top} A \Pi_+$ exist.
            We know from \cite[Section~2]{son2022symplectic} that the negative symplectic eigenvalues of $\Pi_-^{\top} A \Pi_-$ are the negative eigenvalues of $i  \Par{-\Pi_-^{\top} A \Pi_-}^{1/2} J \Par{-\Pi_-^{\top} A \Pi_-}^{1/2}$, and the positive symplectic eigenvalues of $\Pi_+^{\top} A \Pi_+$ are the positive eigenvalues of $i  \Par{\Pi_+^{\top} A \Pi_+}^{1/2} J \Par{\Pi_+^{\top} A \Pi_+}^{1/2}$.
            Therefore, it suffices to show that the non-zero symplectic eigenvalues of $A$ are the non-zero symplectic eigenvalues of $\Pi_-^{\top} A \Pi_-$ and $\Pi_+^{\top} A \Pi_+$ put together.

            Suppose the dimensions of $\operatorname{range}(\Pi_-)$, $\operatorname{range}(\Pi_0)$, $\operatorname{range}(\Pi_+)$ are $2k, 2\ell, 2m$, respectively.
            Thus, $-\Pi_-^{\top} A \Pi_-$ and $\Pi_+^{\top} A \Pi_+$ have ranks $2k$ and $2m$, respectively.
            Let $\mu_1,\ldots, \mu_k$ and $\eta_1,\ldots, \eta_m$ denote the non-zero symplectic eigenvalues of $\Pi_-^{\top} A \Pi_-$ and $\Pi_+^{\top} A \Pi_+$, respectively.
            By Theorem~\ref{thm:williamson_general}, there exist $M, N \in \Sp(2n)$ such that 
            \begin{align}
                M^{\top} \Pi_-^{\top} A \Pi_- M &= D_- \oplus D_-, \label{eq:pi-minus-a-williamson}\\
                N^{\top} \Pi_+^{\top} A \Pi_+ N &= D_+ \oplus D_+.\label{eq:pi-plus-a-williamson}
            \end{align}
            where $D_-$ and $D_+$ are the $n \times n$ diagonal matrices given by $D_-=\diag\Par{\mu_1,\ldots, \mu_k, 0,\ldots, 0}$ and $D_+=\diag\Par{\eta_1,\ldots, \eta_m, 0,\ldots, 0}$.
            Let $w_1,\ldots, w_n, z_1,\ldots, z_n$ be the columns of $M$.
            We have 
            \begin{align}
                \Span \{w_{1},\ldots, w_k, z_{1},\ldots, z_k \}
                    & = \Span \{w_{k+1},\ldots, w_n, z_{k+1},\ldots, z_n\}^{\perp_{\operatorname{s}}} \\
                    & = \Ker\Par{\Pi_-^{\top} A \Pi_-}^{\perp_{\operatorname{s}}} \\
                    & = \Ker\Par{\Pi_-}^{\perp_{\operatorname{s}}} \\
                    & = \ran\Par{\Pi_-}.
            \end{align}
            We know from Proposition~\ref{prop:williamson-alternate-sym-projection} that $\Pi_+ \Pi_-=0$, thus implying for $1\leq i \leq k$ that $w_i, z_i \in \ker(\Pi_+)$.
            Using the fact that $A=\Pi_-^{\top} A \Pi_- + \Pi_+^{\top} A \Pi_+$, we then get
            $\Pi_-^{\top} A \Pi_- w_i = Aw_i$ and $\Pi_-^{\top} A \Pi_- z_i = Az_i$ for all $1\leq i \leq k$.
            The equation \eqref{eq:pi-minus-a-williamson} thus implies for $1\leq i \leq k$:
            \begin{align}
                A w_i &= \mu_i J z_i \\
                A z_i &= - \mu_i J w_i.
            \end{align}
            Let $u_1,\ldots, u_n, v_1,\ldots, v_n$ be the columns of $N$.
            By a similar arguments as given earlier, we get for $j=1,\ldots, m$:
            \begin{align}
                A u_j &= \eta_j J v_j \\
                A v_j &= - \eta_j J u_j.
            \end{align}
            Let $\{x_1,\ldots, x_\ell, y_1,\ldots, y_\ell\}$ be a symplectic basis of $\ran\Par{\Pi_0}$.       
            Let us choose
            \ba{S\coloneqq [w_1,\ldots, w_k, z_1,\ldots, z_k]\diamond [x_1,\ldots, x_\ell, y_1,\ldots, y_\ell] \diamond [u_1,\ldots, u_m, v_1,\ldots, v_m].}
            It is easy to verify that $S \in \Sp(2n)$ and $S^{\top} A S= D \oplus D$, where $D$ is the $n\times n$ diagonal matrix given by $D = \diag\Par{\mu_1,\ldots, \mu_k, 0,\ldots, 0, \eta_1,\ldots, \eta_m}$.
            This completes the proof.
    \end{proof}

\section{Explicit Williamson's decomposition for a subset of $\SpSm(2n)$}\label{sec:special-class}
    For $A \in \operatorname{Sm}(2n)$, let $\mathscr{E}_-, \mathscr{E}_0, \mathscr{E}_+$ denote the eigen subspaces of $A$ spanned by the eigenvectors corresponding to its negative, zero, and positive eigenvalues, respectively.
    We define $\operatorname{EigSpSm}(2n)$ to be the set of those matrices $A \in \operatorname{Sm}(2n)$ for which $\mathscr{E}_-, \mathscr{E}_0, \mathscr{E}_+$ are pairwise symplectically orthogonal symplectic subspaces, and each of these subspaces is invariant under $JA$.
    % Observe that $\Pd(2n) \subset \operatorname{SpPsd}(2n) \subset \operatorname{EigSpSm}(2n) \subset \SpSm(2n)$.

    In this section, we provide an explicit description of symplectic eigenvalues and diagonalizing symplectic matrices in Williamson's decomposition for matrices in $\operatorname{EigSpSm}(2n)$.
    Furthermore, we establish perturbation bounds on the symplectic eigenvalues of matrices in $\operatorname{EigSpSm}(2n)$.

    We begin with some preliminary results that will be helpful in the subsequent parts of the section.
    \begin{boxed}
    \begin{lemma}\label{lem:positive_definite_restriction}
    Let $A \in \operatorname{S}(2n)$ and $\mathscr{E}$ be a symplectic subspace of $\mathds{R}^{2n}$ of dimension $2 k$.
    Suppose $\mathscr{E}$ is an invariant subspace of both $JA$ and $A$, and that $A$ is positive definite on $\mathscr{E}$.
    Then there exist $k$ positive numbers $\gamma_1,\ldots, \gamma_k$ and a symplectic basis $\{u_1,\ldots, u_k, v_1,\ldots, v_k\}$ of $\mathscr{E}$ such that for all $1 \leq i \leq k$,
    \begin{align}
        A u_i &= \gamma_i J v_i, \\
        A v_i &= -\gamma_i J u_i.
    \end{align}
    \end{lemma}
    \end{boxed}
    \begin{proof}
    %We shall extend $\restr{A}{\mathscr{E}}$ to a positive definite matrix $\widetilde{A} \in \mathds{P}_{2n}(\mathds{R})$.  
    Let $P$ be the orthogonal projection onto the subspace $\mathscr{E}$.
    Set $\hat{A}\coloneqq A+I-P$, where $I$ is the identity matrix.
    Let $x \in \mathscr{E}$ and $x^{\perp} \in \mathscr{E}^{\perp}$ be arbitrary.
    We have 
    \begin{align}
        \hat{A}(x+x^{\perp}) 
            &=Ax+x-Px + Ax^{\perp}+x^{\perp}-Px^{\perp} \\
            &= Ax+x-x + x^{\perp} \\
            &= Ax+x^{\perp}.
    \end{align}
    If $x + x^{\perp} \neq 0$, i.e., $x \neq 0$ or $x^{\perp} \neq 0$.
    This then implies
    \begin{align}
        \langle x+x^{\perp}, \hat{A}(x+x^{\perp}) \rangle
            &= \langle x+x^{\perp}, Ax+x^{\perp} \rangle \\
            &= \langle x, Ax \rangle + \langle x, x^{\perp} \rangle + \langle x^{\perp}, Ax \rangle + \langle x^{\perp}, x^{\perp} \rangle \\
            &= \langle x, Ax \rangle + \langle x^{\perp}, x^{\perp} \rangle > 0.
    \end{align}
    This implies that $\hat{A} \in \operatorname{Pd}(2n)$, and it is easy to see that $\mathscr{E}$ is invariant under $J\hat{A}$.
    By Proposition~4.1 of \cite{mishra2024equality}, there exists a symplectic basis $\{u_1,\ldots, u_k, v_1,\ldots, v_k\}$ of $\mathscr{E}$ such that
    \begin{align}
        \hat{A} u_i &= \gamma_i J v_i, \\
        \hat{A} v_i &= -\gamma_i J u_i.
    \end{align}
    By definition, we have $\hat{A} u_i= Au_i$ and $\hat{A} v_i= Av_i$ for all $1\leq i \leq k$.
    This completes the proof.
    % Define $\widetilde{A} \in \operatorname{M}_{2n}(\mathds{R})$ whose columns are given by
    % \begin{align}\label{eq:columns_of_extension}
    %     \widetilde{A}_i \coloneqq \left[AP \right]_i+ \left[I-P \right]_i, \quad \text{for } i=1,\ldots, n.
    % \end{align}
    % The matrix $\widetilde{A}$ is the canonical matrix of the linear transformation $T : \mathds{R}^{2n}\to \mathds{R}^{2n}$ defined as follows: for $x \in \mathscr{E}$ and $y \in \mathscr{E}^{\perp}$,
    % \begin{align}
    %     T(x \oplus y) \coloneqq Ax + y.
    % \end{align}
    \end{proof}
    The following well-known result on commuting normal matrices plays key role in constructing symplectic matrices in Williamson's decomposition for matrices in $\operatorname{EigSpSm}(2n)$.
    See Theorem~2.5.15 of \cite{horn2012matrix} for a proof.
    \begin{boxed}
        \begin{lemma}\label{lem:commuting-matrices-real}
            Let $A, B \in \operatorname{M}(n)$ be normal matrices.
            If $A$ and $B$ commute, then there exists $P \in \operatorname{O}(n)$ and a non-negative integer $r$ such that $P^{\top} A P$ and $P^{\top} B P$ are block-diagonal matrices of the form:
            \begin{align}
                P^{\top} AP 
                    &=  
                     \Delta_1 \oplus \begin{pmatrix} \alpha_1 & \beta_1 \\ -\beta_1 & \alpha_1\end{pmatrix} \oplus \cdots \oplus \begin{pmatrix} \alpha_r & \beta_r \\ -\beta_r & \alpha_r\end{pmatrix}   , \\
                P^{\top}BP
                    &=  
                     \Delta_2 \oplus \begin{pmatrix} \gamma_1 & \delta_1 \\ -\delta_1 & \gamma_1\end{pmatrix} \oplus \cdots \oplus \begin{pmatrix} \gamma_r & \delta_r \\ -\delta_r & \gamma_r\end{pmatrix},
            \end{align}
            where $\Delta_1, \Delta_2 \in \operatorname{M}(n-2r)$ are diagonal matrices; $\alpha_i, \beta_i, \gamma_i, \delta_i$ are real numbers for all $1 \leq i \leq r$; and for each $i\in \{1,\ldots,r\}$, $\beta_i>0$ or $\delta_i>0$.
        \end{lemma}
    \end{boxed}

\subsection{Description of symplectic eigenvalues for $\operatorname{EigSpSm}(2n)$}
    \label{sec:williamson-construction-1}
    
    The symplectic eigenvalues of a matrix $A \in \operatorname{EigSpSm}(2n)$ are given by a combination of negative and non-negative eigenvalues of the Hermitian matrices $i  A_-^{1/2}JA^{1/2}_-$ and $i  A_+^{1/2}JA^{1/2}_+$ as stated below.
    \begin{boxed}
    \begin{theorem}\label{thm:williamson_gen_special}
    The symplectic eigenvalues of $A \in \operatorname{EigSpS}(2n)$ are given by $\frac{1}{2} \xi(A)$ zeros, the negative eigenvalues of $i  A_-^{1/2}JA^{1/2}_-$, and the positive eigenvalues of $i  A_+^{1/2}JA^{1/2}_+$.
    \end{theorem}
    \end{boxed}
    \begin{proof}
    Let $\mathscr{E}_-, \mathscr{E}_0, \mathscr{E}_+$ denote the eigen subspaces of $A$ spanned by the eigenvectors corresponding to its negative, zero, and positive eigenvalues, respectively. 
    Let $\Pi_{-}, \Pi_0, \Pi_+$ denote the orthogonal projections onto $\mathscr{E}_-, \mathscr{E}_0, \mathscr{E}_+$, respectively.
    By definition, $\Pi_{-}, \Pi_0, \Pi_+$ are also symplectic orthogonal projections onto the symplectic subspaces $\mathscr{E}_-, \mathscr{E}_0, \mathscr{E}_+$, respectively.
    Also, we have
    \begin{align}
        \Pi_{-}^{\top} A \Pi_- &= -A_-, \\
        \Pi_{+}^{\top} A \Pi_+ &= A_+.
    \end{align}
    By Proposition~\ref{prop:sym-eigenvalues-eigenvalues-relation}, the negative eigenvalues of $i  A_-^{1/2}JA^{1/2}_-$ and the positive eigenvalues of $i  A_+^{1/2}JA^{1/2}_+$, together with $\frac{1}{2} \xi(A)$ zeros are the symplectic eigenvalues of $A$.
    \end{proof}

    \subsection{Description of symplectic matrices in Williamson's decomposition for $\operatorname{EigSpSm}(2n)$}
    \label{sec:williamson-construction-2}
    Let $A \in \operatorname{EigSpSm}(2n)$.
    In what follows, we explicitly construct a symplectic matrix that diagonalizes $A$ in the sense of Williamson's theorem.
    
    We know that the matrices $A_{-}^{1/2}$ and $A_{+}^{1/2}$ commute with each other and satisfy $A_{-}^{1/2}A_{+}^{1/2}=0$.
    Therefore, the skew-symmetric matrices $A_{-}^{1/2}JA_{-}^{1/2}$ and $A_{+}^{1/2}JA_{+}^{1/2}$ commute with each other and their product is equal to zero.
    By Lemma~\ref{lem:commuting-matrices-real}, there exists $U \in \operatorname{O}(2n)$ and a non-negative integer $r$ such that
    \begin{align}
        U^{\top}A_{-}^{1/2}JA_{-}^{1/2}U 
            &=  
            \Delta_1 \oplus \begin{pmatrix} \alpha_1 & \beta_1 \\ -\beta_1 & \alpha_1\end{pmatrix} \oplus \cdots \oplus \begin{pmatrix} \alpha_r & \beta_r \\ -\beta_r & \alpha_r\end{pmatrix}  , \label{eq:j-aplus-j-real-schur-form} \\
        U^{\top}A_{+}^{1/2}JA_{+}^{1/2}U
            &=  
            \Delta_2 \oplus \begin{pmatrix} \gamma_1 & \delta_1 \\ -\delta_1 & \gamma_1\end{pmatrix} \oplus \cdots \oplus \begin{pmatrix} \gamma_r & \delta_r \\ -\delta_r & \gamma_r\end{pmatrix},\label{eq:j-aminus-j-real-schur-form}
    \end{align}
    where $\Delta_1, \Delta_2$ are real diagonal matrices of size $(2n-2r) \times (2n-2r)$; the parameters $\alpha_i, \beta_i, \gamma_i, \delta_i$ are real numbers such that $\beta_i>0$ or $\delta_i>0$ for all $1\leq i \leq r$.
    Since both $U^{\top}A_{-}^{1/2}JA_{-}^{1/2}U$ and $U^{\top}A_{+}^{1/2}JA_{+}^{1/2}U$ are real skew-symmetric matrices, their diagonal elements are zero whence $\Delta_1=\Delta_2=0$ and $\alpha_i=\gamma_i=0$ for all $1\leq i \leq r$.
    The fact that the product of the matrices in the left-hand sides of \eqref{eq:j-aplus-j-real-schur-form} and \eqref{eq:j-aminus-j-real-schur-form} is zero implies that $\beta_i \delta_i=0$.
    This implies that for all $1\leq i \leq r$, exactly one of $\beta_i$ and $\delta_i$ is positive.

    We know that the kernel of $A_-$ is $\mathscr{E}_0 + \mathscr{E}_+$, which is a symplectic subspace of $\mathds{R}^{2n}$ of dimension $2(\ell+m)$.
    It is shown in \cite[Section~2]{son2022symplectic} that $\operatorname{ker}(A_{-}^{1/2}JA_{-}^{1/2})=\operatorname{ker}(A_-)$, which implies $\operatorname{rank}(A_{-}^{1/2}JA_{-}^{1/2})=2k$.
    Similarly, we get $\operatorname{rank}(A_{+}^{1/2}JA_{+}^{1/2})=2m$.
    Therefore, we must have $r=k+m$, there exist distinct indices $1 \leq i_1< \cdots < i_k \leq k+m$ and $1 \leq j_1< \cdots < j_m \leq k+m$ such that for $i \in \{ i_1,\ldots, i_k\}$, we have $\beta_i>0, \delta_i =0$ and for $j \in \{ j_1,\ldots, j_m \}$, we have $\beta_j=0, \delta_j >0$.
    Let $D_-$ and $D_+$ be $2n \times 2n$ diagonal matrices whose $i$th diagonal entries are given by
    \begin{align}
        (D_-)_i 
            &=
            \begin{cases}
                0 & \text{if } i \in \{1,\ldots,\ell\} \cup \{\ell+j_1,\ldots, \ell+j_m\}, \\
                \beta_{i-\ell} & \text{if } i \in \{ \ell+i_1,\ldots, \ell+i_k\}, \\
            \end{cases} \\
        (D_+)_i 
            &=
            \begin{cases}
                0 & \text{if } i \in \{1,\ldots,\ell\} \cup \{\ell+i_1,\ldots, \ell+i_k\}, \\
                \delta_{i-\ell} & \text{if } i \in \{\ell+ j_1,\ldots, \ell+ j_m\}.
            \end{cases}
    \end{align}
    Let $e_1,\ldots, e_{2n}$ denote the standard unit vectors in $\mathds{R}^{2n}$.
    Let $P$ denote the permutation matrix $[e_1, e_{3}, \ldots, e_{2n-1}, e_{2}, e_4,\ldots, e_{2n}]$.
    We then get
    \begin{align}
        P^{\top}U^{\top}A_{-}^{1/2}JA_{-}^{1/2}UP
            &=  (D_- \oplus D_-) J , \label{eq:j-aplus-j-real-schur-form-simplification} \\
        P^{\top}U^{\top}A_{+}^{1/2}JA_{+}^{1/2}UP
            &=  
            (D_+ \oplus D_+) J.\label{eq:j-aminus-j-real-schur-form-simplification}
    \end{align}
    Let $\Pi_-$ and $\Pi_+$ denote the following isometries
    \begin{align}
        \Pi_- &\coloneqq [e_{\ell+i_1}, \ldots, e_{\ell+i_k}, e_{n+\ell+i_1}, \ldots, e_{n+\ell+i_k}], \\
        \Pi_+ &\coloneqq [e_{\ell+j_1}, \ldots, e_{\ell+j_m}, e_{n+\ell+j_1}, \ldots, e_{n+\ell+j_m}].
    \end{align}
    From \eqref{eq:j-aplus-j-real-schur-form-simplification} and \eqref{eq:j-aminus-j-real-schur-form-simplification} we thus get
    \begin{align}
        \Pi_-^{\top}P^{\top}U^{\top}A_{-}^{1/2}J_{2n}A_{-}^{1/2}UP\Pi_-
            &=  (\widetilde{D}_- \oplus \widetilde{D}_-) J_{2k} , \label{eq:j-aplus-j-real-schur-form-simplification-2} \\
        \Pi_+^{\top}P^{\top}U^{\top}A_{+}^{1/2}J_{2n}A_{+}^{1/2}UP\Pi_+
            &=  
            (\widetilde{D}_+ \oplus \widetilde{D}_+) J_{2m},\label{eq:j-aminus-j-real-schur-form-simplification-2}
    \end{align}
    where $\widetilde{D}_- \coloneqq \operatorname{diag}(\beta_{i_1},\ldots, \beta_{i_k})$ and $\widetilde{D}_+ \coloneqq \operatorname{diag}(\delta_{j_1},\ldots, \delta_{j_m})$.
    Choose 
    \begin{align}
        \widetilde{M}_-
            &\coloneqq J_{2n}A_{-}^{1/2}UP\Pi_- \Par{\widetilde{D}_-^{-1/2} \oplus \widetilde{D}_-^{-1/2}} J_{2k}^{\top}, \label{eq:def_sym_matrix_aplus}\\
        \widetilde{M}_+
            &\coloneqq J_{2n}A_{+}^{1/2}UP\Pi_+ \Par{\widetilde{D}_+^{-1/2} \oplus \widetilde{D}_+^{-1/2}} J_{2m}^{\top}. \label{eq:def_sym_matrix_aminus}
    \end{align}
    It is easy to see from \eqref{eq:j-aplus-j-real-schur-form-simplification-2} and \eqref{eq:j-aminus-j-real-schur-form-simplification-2} that $\widetilde{M}_- \in \operatorname{Sp}(2n, 2k)$ and $\widetilde{M}_+ \in \operatorname{Sp}(2n, 2m)$.
    We observe that $A_{-}^{1/2}JA=-A_{-}^{1/2}JA_-$, which follows from the fact that $\mathscr{E}_-$ and $\mathscr{E}_+$ are invariant under $JA$.
    Therefore, we get
    \begin{align}
        &\widetilde{M}_-^{\top} A \widetilde{M}_- \\
            &=-J_{2k} \Par{\widetilde{D}_-^{-1/2} \oplus \widetilde{D}_-^{-1/2}} \Pi_-^{\top}P^{\top}U^{\top}A_{-}^{1/2}  J_{2n}A J_{2n}A_{-}^{1/2}UP\Pi_- \Par{\widetilde{D}_-^{-1/2} \oplus \widetilde{D}_-^{-1/2}} J_{2k}^{\top} \\
            &=J_{2k} \Par{\widetilde{D}_-^{-1/2} \oplus \widetilde{D}_-^{-1/2}} \Pi_-^{\top}P^{\top}U^{\top}A_{-}^{1/2}  J_{2n}A_- J_{2n}A_{-}^{1/2}UP\Pi_- \Par{\widetilde{D}_-^{-1/2} \oplus \widetilde{D}_-^{-1/2}} J_{2k}^{\top}\\
            &=J_{2k} \Par{\widetilde{D}_-^{-1/2} \oplus \widetilde{D}_-^{-1/2}} \Pi_-^{\top}\left(P^{\top}U^{\top} A_{-}^{1/2}  J_{2n} A_{-}^{1/2} UP \right)^2\Pi_- \Par{\widetilde{D}_-^{-1/2} \oplus \widetilde{D}_-^{-1/2}} J_{2k}^{\top} \\
            &=-J_{2k} \Par{\widetilde{D}_-^{-1/2} \oplus \widetilde{D}_-^{-1/2}} \Pi_-^{\top} (D_-^2 \oplus D_-^2) \Pi_- \Par{\widetilde{D}_-^{-1/2} \oplus \widetilde{D}_-^{-1/2}} J_{2k}^{\top} \\
            &=-J_{2k} \Par{\widetilde{D}_-^{-1/2} \oplus \widetilde{D}_-^{-1/2}}  (\widetilde{D}_-^2 \oplus \widetilde{D}_-^2) \Par{\widetilde{D}_-^{-1/2} \oplus \widetilde{D}_-^{-1/2}} J_{2k}^{\top} \\
            &=- J_{2k}(\widetilde{D}_- \oplus \widetilde{D}_-) J_{2k}^{\top} \\
            &= - (\widetilde{D}_- \oplus \widetilde{D}_-).\label{eq:sym_diagonalization_a_minus}
    \end{align}
     By similar arguments, one can show that
     \begin{align}
         \widetilde{M}_+^{\top} A \widetilde{M}_+
            &= \widetilde{D}_+ \oplus \widetilde{D}_+.\label{eq:sym_diagonalization_a_plus}
     \end{align}
    Choose any $\widetilde{M}_0 \in \operatorname{Sp}(2n, 2\ell)$ whose columns form a symplectic basis of $\mathscr{E}_0$.
    Define
    \begin{align}
        \widetilde{M} \coloneqq \widetilde{M}_0 \diamond \widetilde{M}_- \diamond \widetilde{M}_+.
    \end{align}
    The matrix $\widetilde{M}$ is symplectic.
    Indeed, we have $A_{+}^{1/2}JA_{-}^{1/2}=0$ since $\mathscr{E}_-$ and $\mathscr{E}_+$ are invariant under $JA$.
    We thus get from \eqref{eq:def_sym_matrix_aplus} and \eqref{eq:def_sym_matrix_aminus} that
    \begin{align}
        \widetilde{M}_+^{\top} J_{2n} \widetilde{M}_- &= 0_{2m, 2k}.\label{eq:sym_condition_tilde_M_1}
    \end{align}
    Since the subspaces $\mathscr{E}_0$ and $\mathscr{E}_-$ are perpendicular to each other, and $\operatorname{range}(A_{-}^{1/2})= \operatorname{range}(A_-) = \mathscr{E}_-$, we get 
    \begin{align}
        \widetilde{M}_0^{\top} J_{2n} \widetilde{M}_- &=0_{2\ell, 2k}.\label{eq:sym_condition_tilde_M_2}
    \end{align}
    By similar arguments, we also get
    \begin{align}
        \widetilde{M}_0^{\top} J_{2n} \widetilde{M}_+ &=0_{2\ell, 2m}.\label{eq:sym_condition_tilde_M_3}
    \end{align}
    The conditions \eqref{eq:sym_condition_tilde_M_1}, \eqref{eq:sym_condition_tilde_M_2}, \eqref{eq:sym_condition_tilde_M_3} thus imply that $\widetilde{M} \in \operatorname{Sp}(2n)$.
    See \cite[Section~2.3]{mishra2023}.
    By \eqref{eq:sym_diagonalization_a_minus}, \eqref{eq:sym_diagonalization_a_plus}, and the fact that $A M_0=0_{2n, 2\ell}$, we get
    \begin{align}
        \widetilde{M}^{\top} A \widetilde{M} &= D \oplus D,
    \end{align}
    where $D\coloneqq (-\widetilde{D}_-) \oplus 0_{\ell, \ell} \oplus \widetilde{D}_+$.
\subsection{Perturbation bounds on symplectic eigenvalues for $\operatorname{EigSpSm}(2n)$}\label{sec:perturbation_bounds}
    In this subsection, we provide perturbation bounds on symplectic eigenvalues of matrices in $\operatorname{EigSpSm}(2n)$ given by Theorem~\ref{thm:williamson_gen_special}.
    These perturbation bounds generalize the known perturbation bounds on symplectic eigenvalues of positive definite matrices given in \cite{bhatia2015symplectic}.
    
    Let $\operatorname{M}(n, \mathds{C})$ denote the set of $n \times n$ complex matrices, and $\operatorname{U}(n, \mathds{C})$ denote the set of $n \times n$ complex unitary matrices.
    A norm $|\!|\!| \cdot |\!|\!|$ on $\operatorname{M}(n, \mathds{C})$ is called unitarily invariant if $|\!|\!|U X V |\!|\!|=|\!|\!| X|\!|\!|$ for all $X \in \operatorname{M}(n, \mathds{C})$ and $U, V \in \operatorname{U}(n, \mathds{C})$.
    For $X, Y , Z \in \operatorname{M}(n, \mathds{C})$, every unitarily invariant norm satisfies $|\!|\!| X Y Z |\!|\!| \leq \|X \| \cdot |\!|\!| Y |\!|\!| \cdot \|Z\|$.
    Here $\|\cdot \|$ denotes the matrix operator norm.
    See Proposition~IV.2.4 of \cite{ma_bhatia}.
    For $A, B \in \Psd(n)$, the following inequality holds \cite[Theorem~X.1.3]{ma_bhatia}:
    \begin{align}
        |\!|\!| A^{1/2}-B^{1/2} |\!|\!| \leq  |\!|\!| |A-B|^{1/2} |\!|\!|.\label{eq:abs_uni_inv_norm_pert_bound_sqrt}
    \end{align}
    % For every $X, Y \in \operatorname{M}(n)$, the following inequality holds (see \cite[Theorem~X.2.1]{ma_bhatia}):
    % \begin{align}
    %     |\!|\!| |A|-|B| |\!|\!| \leq \sqrt{2} \left(|\!|\!| A+B |\!|\!| \cdot |\!|\!| A-B |\!|\!| \right)^{1/2}.\label{eq:abs_uni_inv_norm_pert_bound}
    % \end{align}

    Given an $n \times n$ complex Hermitian matrix $X$, let $\lambda^{}(X)$ denote the $n$-vector consisting of the eigenvalues of $X$ arranged in the \emph{decreasing}  order.  
    Let $\operatorname{Eig}(X)$ denote the $n \times n$ diagonal matrix whose diagonal elements are given by the entries of $\lambda^{}(X)$. 
    The Lidskii--Wielandt theorem \cite[IV.62]{ma_bhatia} gives
    \begin{align}
        |\!|\!| \operatorname{Eig}(X)-\operatorname{Eig}(Y) |\!|\!| \leq |\!|\!| X-Y |\!|\!|.\label{eq:lidskii-wielandt}
    \end{align}
    For $A \in \operatorname{S}(2n)$, let $\widehat{D}(A)$ be the $2n \times 2n$ diagonal matrix
    \begin{align}
        \widehat{D}(A) \coloneqq \operatorname{Eig}^{}(|A_+^{1/2}J_{2n} A_+^{1/2}|) + \operatorname{Eig}^{}(-|A_-^{1/2}J_{2n} A_-^{1/2}|).\label{eq:diagonal-d-tilde}
    \end{align}
    Since the eigenvalues of $i  A_+^{1/2}J_{2n} A_+^{1/2}$ and $i  A_-^{1/2}J_{2n} A_-^{1/2}$ occur in pairs of negative-positive, the diagonal elements of $\widehat{D}(A)$ occur in pairs of equal entries, and we denote the diagonal elements of $\widehat{D}(A)$ by $d_1(A), d_1(A), \ldots, d_n(A), d_n(A)$.
    
     The next lemma gives a perturbation bound on $\widehat{D}(A)$.    
    We know from Theorem~\ref{thm:williamson_gen_special} that if $A \in \operatorname{EigSpSm}(2n)$, then the diagonal elements of $\widehat{D}(A)$ are the symplectic eigenvalues of $A$ given by Theorem~\ref{thm:williamson_gen_special}, each counted twice.
    \begin{boxed}
    \begin{proposition}\label{prop:perturbation_bounds}
        Let $A, B \in \operatorname{EigSpSm}(2n)$.
        We have
        \begin{multline}
            |\!|\!| \widehat{D}(A) - \widehat{D}(B) |\!|\!| 
                \leq \left(\|A_+^{1/2} \| + \| B_+^{1/2}\| \right) |\!|\!|  |A_+  -  B_+|^{1/2} |\!|\!| \\ + 
                \left(\|A_-^{1/2} \| + \| B_-^{1/2}\| \right) |\!|\!|  |A_-  -  B_-|^{1/2} |\!|\!|.\label{eq:perturbation_bound_main}
        \end{multline}
        In the special cases of the operator norm and the Frobenius norm, we get
        \begin{align}
            \max_{1\leq i \leq n} |d_i(A)-d_i(B)| 
                &\leq \left(\|A_+^{1/2} \| + \| B_+^{1/2}\| \right) \|  A_+  -  B_+\|^{1/2} \notag \\
                &\hspace{0.5cm}+ \left(\|A_-^{1/2} \| + \| B_-^{1/2}\| \right) \|  A_-  -  B_-\|^{1/2},\label{eq:perturbation_bound_operatornorm} \\
            \sqrt{2}\left(\sum_{ i =1}^n |d_i(A)-d_i(B)|^2 \right)^{1/2} 
                &\leq \left(\|A_+^{1/2} \| + \| B_+^{1/2}\| \right) \operatorname{Tr}(  |A_+  -  B_+|)^{1/2} \notag \\ 
                &\hspace{0.5cm}+ \left(\|A_-^{1/2} \| + \| B_-^{1/2}\| \right) \operatorname{Tr} ( |A_-  -  B_-|)^{1/2}.\label{eq:perturbation_bound_frobeniusnorm}
        \end{align}
    \end{proposition}
    \end{boxed}
    \begin{proof}
        By definition \eqref{eq:diagonal-d-tilde} and triangle inequality, we get
        \begin{multline}
            |\!|\!| \widehat{D}(A) - \widehat{D}(B) |\!|\!| 
            \leq |\!|\!| \operatorname{Eig}^{}(|A_+^{1/2}J_{2n} A_+^{1/2}|)  - \operatorname{Eig}^{}(|B_+^{1/2}J_{2n} B_+^{1/2}|)|\!|\!| \\
             + |\!|\!|  \operatorname{Eig}^{}(-|A_-^{1/2}J_{2n} A_-^{1/2}|)  - \operatorname{Eig}^{}(-|B_-^{1/2}J_{2n} B_-^{1/2}|) |\!|\!|.\label{eq:sp_eig_perturbation_inequality_1}
        \end{multline}
        We know that the eigenvalues of $i  A_+^{1/2}J_{2n} A_+^{1/2}$ and $i  B_+^{1/2}J_{2n} B_+^{1/2}$ occur in positive negative pairs. 
        Therefore, using the unitary invariance of the norm, we get
        \begin{multline}
            |\!|\!| \operatorname{Eig}^{}(|A_+^{1/2}J_{2n} A_+^{1/2}|)  - \operatorname{Eig}^{}(|B_+^{1/2}J_{2n} B_+^{1/2}|)|\!|\!| \\ = |\!|\!| \operatorname{Eig}^{}(i  A_+^{1/2}J_{2n} A_+^{1/2})  - \operatorname{Eig}^{}(i  B_+^{1/2}J_{2n} B_+^{1/2})|\!|\!|.\label{eq:norm_equality_by_unitary_invariance_1}
        \end{multline}
        Similarly, we also have
        \begin{multline}
            |\!|\!|  \operatorname{Eig}^{}(-|A_-^{1/2}J_{2n} A_-^{1/2}|)  - \operatorname{Eig}^{}(-|B_-^{1/2}J_{2n} B_-^{1/2}|) |\!|\!| \\ = |\!|\!| \operatorname{Eig}^{}(i  A_-^{1/2}J_{2n} A_-^{1/2})  - \operatorname{Eig}^{}(i  B_-^{1/2}J_{2n} B_-^{1/2}) |\!|\!|.\label{eq:norm_equality_by_unitary_invariance_2}
        \end{multline}
        Substituting \eqref{eq:norm_equality_by_unitary_invariance_1} and \eqref{eq:norm_equality_by_unitary_invariance_2} into the right-hand side of \eqref{eq:sp_eig_perturbation_inequality_1}, we get
        \begin{multline}
            |\!|\!| \widehat{D}(A) - \widehat{D}(B) |\!|\!| 
            \leq |\!|\!| \operatorname{Eig}^{}(i  A_+^{1/2}J_{2n} A_+^{1/2})  - \operatorname{Eig}^{}(i  B_+^{1/2}J_{2n} B_+^{1/2})|\!|\!| \\
             + |\!|\!|  \operatorname{Eig}^{}(i  A_-^{1/2}J_{2n} A_-^{1/2})  - \operatorname{Eig}^{}(i  B_-^{1/2}J_{2n} B_-^{1/2}) |\!|\!|.\label{eq:sp_eig_perturbation_inequality_2}
        \end{multline}
        We now apply the same arguments as given in the proof of Theorem~7 of \cite{bhatia2015symplectic} to bound each term in the right-hand side of \eqref{eq:sp_eig_perturbation_inequality_2}.
        
        By the Lidskii--Wielandt theorem \eqref{eq:lidskii-wielandt} and the relation \eqref{eq:abs_uni_inv_norm_pert_bound_sqrt}, we get
        \begin{align}
            &|\!|\!| \operatorname{Eig}^{}(i  A_+^{1/2}J_{2n} A_+^{1/2})  - \operatorname{Eig}^{}(i  B_+^{1/2}J_{2n} B_+^{1/2})|\!|\!| \notag \\
            &\leq |\!|\!|  A_+^{1/2}J_{2n} A_+^{1/2}  -  B_+^{1/2}J_{2n} B_+^{1/2}|\!|\!| \\
            &\leq |\!|\!|  A_+^{1/2}J_{2n} A_+^{1/2}  -  A_+^{1/2}J_{2n} B_+^{1/2}|\!|\!| + |\!|\!|  A_+^{1/2}J_{2n} B_+^{1/2}  -  B_+^{1/2}J_{2n} B_+^{1/2}|\!|\!| \\
            &= |\!|\!|  A_+^{1/2}J_{2n} (A_+^{1/2}  -  B_+^{1/2} )|\!|\!| + |\!|\!|  (A_+^{1/2}  -  B_+^{1/2} )J_{2n} B_+^{1/2}|\!|\!| \\
            &\leq \|A_+^{1/2}J_{2n} \|\cdot |\!|\!|  A_+^{1/2}  -  B_+^{1/2} |\!|\!| + |\!|\!|  A_+^{1/2}  -  B_+^{1/2} |\!|\!|  \cdot \|J_{2n} B_+^{1/2}\| \\
            &= \left(\|A_+^{1/2} \| + \| B_+^{1/2}\| \right) |\!|\!|  A_+^{1/2}  -  B_+^{1/2} |\!|\!| \\
            &\leq \left(\|A_+^{1/2} \| + \| B_+^{1/2}\| \right) |\!|\!|  |A_+  -  B_+|^{1/2} |\!|\!|.\label{eq:inequality-positive-part}
        \end{align}
        Similarly,
        \begin{align}
            |\!|\!| \operatorname{Eig}^{}(i  A_-^{1/2}J_{2n} A_-^{1/2})  - \operatorname{Eig}^{}(i  B_-^{1/2}J_{2n} B_-^{1/2})|\!|\!| 
            \leq \left(\|A_-^{1/2} \| + \| B_-^{1/2}\| \right) |\!|\!|  |A_-  -  B_-|^{1/2} |\!|\!|.\label{eq:inequality-negative-part}
        \end{align}
        Substituting \eqref{eq:inequality-positive-part} and \eqref{eq:inequality-negative-part} into \eqref{eq:sp_eig_perturbation_inequality_2} gives the desired perturbation bound \eqref{eq:perturbation_bound_main}.
        The other perturbation bounds \eqref{eq:perturbation_bound_operatornorm} and \eqref{eq:perturbation_bound_frobeniusnorm} follow directly from \eqref{eq:perturbation_bound_main}.
    \end{proof}

    \begin{boxed}
    \begin{remark}\emph{
        In Proposition~\ref{prop:perturbation_bounds}, if the matrices $A$ and $B$ are positive definite, then we have $A_-=B_-=0$, $A_+=A$, and $B_+=B$. 
        The perturbation bound \eqref{eq:perturbation_bound_main} in this case reduces to the perturbation bound of symplectic eigenvalues of $A$ and $B$ given in Theorem~7 of \cite{bhatia2015symplectic}.
    }\end{remark}
    \end{boxed}

 \section{Interpretations of symplectic orthogonal projection and some of the results for quadratic forms on general symplectic spaces}
 \label{sec:interpretation}
    We first recall some basic theory of quadratic forms and symplectic geometry.
    \vspace{0.5cm}
    
    \noindent \textbf{\textit{Quadratic forms}.}
    A quadratic form on a real vector space $\mathscr{X}$ is a map $Q:\mathscr{X}\to\mathds{R}$
    that satisfies $(i)$ \textit{Homogeneity of order two:} $Q(c x)=c^2 Q(x)$ for $c \in \mathds{R}$ and $x \in \mathscr{X}$, and $(ii)$ \textit{Polar identity:} the map $(x,y) \mapsto \Phi_Q(x, y) \coloneqq \frac{1}{2}\left(Q(x+y)-Q(x)-Q(y)\right)$ is a symmetric bilinear form.
    It is straightforward to verify that the mapping $Q \mapsto \Phi_Q$ is a one-to-one correspondence between the set of quadratic forms and the set of symmetric bilinear forms on $\mathscr{X}$.
    %Equivalently, $Q$ is a quadratic form on $\mathscr{V}$ if, in a fixed coordinate system, $Q(v)$ is a homogeneous polynomial of degree two of the coordinates of $v$.
    % $Q$ is said to be positive semi-definite if $\Phi_Q$ is positive semi-definite.
    % Similarly, $Q$ is said to be positive definite if $\Phi_Q$ is positive definite.
    % See the first two chapters of \cite{elman2008algebraic} for a detailed theory of bilinear and quadratic forms.
    If $\mathscr{X}$ is $n$-dimensional, then $\Phi_Q$ can be represented by an $n \times n$ symmetric matrix in a given basis of $\mathscr{V}$.
    By Sylvester's law of inertia, the inertia of any symmetric matrix representing $\Phi_Q$ is independent of the choice of the basis.
    We denote by $\nu(Q), \xi(Q), \pi(Q)$, respectively, the number of positive, zero, and negative eigenvalues of a symmetric matrix representing the bilinear form $\Phi_Q$ in a basis.
    
    \vspace{0.5cm}
    \noindent \textbf{\textit{Hamiltonian map and complex structure}.}
    Let $\left(\mathscr{V}, \omega \right)$ be a symplectic space.
    Associated with every quadratic form $Q$ on $\mathscr{V}$ is a unique linear map $H_Q: \mathscr{V} \to \mathscr{V}$ given by
    \begin{align}
    \label{eq:bilinear_Hamilton_relation}
        \Phi_Q(u, v) = \omega(u, H_Q(v)), \qquad  u,v \in \mathscr{V}.
    \end{align}
    The map $H_Q$ is known as the \emph{Hamilton map} of $Q$ (see, e.g.,  \cite{OTTOBRE20124000}).
    There exists an automorphism $J: \mathscr{V} \to \mathscr{V}$, called a \emph{complex structure compatible with $\omega$} \cite[Lemma~2.5.5]{Mcduff_salamon}, satisfying the following conditions
    \begin{itemize}
        \item $J^2=-\mathds{1}$, where $\mathds{1}$ is the identity map,
        \item $\omega(Ju, Jv)=\omega(u,v)$ for all $u, v \in \mathscr{V}$, and
        \item $g_J(u, v) \coloneqq \omega(u, Jv)$ defines an inner product on $\mathscr{V}$.
    \end{itemize}
    The space of complex structures can be identified with the Siegel upper half space \cite[Lemma~2.5.12]{Mcduff_salamon}.
    It is known that there exists a symplectic basis of $\left(\mathscr{V}, \omega \right)$ which is also an orthonormal basis of the inner product space $\left(\mathscr{V}, g_J \right)$ \cite[Lemma~2.4.5]{Mcduff_salamon}.
    We call such a basis \emph{$J$-orthosymplectic basis} of $\left(\mathscr{V}, \omega \right)$.
    The standard basis of $\mathds{R}^{2n}$ is an example of a $J$-orthosymplectic basis for $J=\begin{psmallmatrix}
             0 & I \\
             -I & 0
        \end{psmallmatrix}.$

    \vspace{0.5cm}
    \noindent \textbf{\textit{Symplectic orthogonal complement}.}
    Let $\mathscr{W}$ be a linear subspace of the symplectic space $\left(\mathscr{V}, \omega \right)$.
    The \emph{symplectic orthogonal complement} of $\mathscr{W}$ is defined as
    \begin{align}
        \mathscr{W}^{\perp_{\operatorname{s}}} \coloneqq \big\{v \in \mathscr{V}: \omega(v, w)=0 \ \forall w \in \mathscr{W} \big\}.
    \end{align}
    Moreover, $\mathscr{W}^{\perp_{\operatorname{s}}}$ is also a linear subspace, and satisfies
    \begin{align}\label{eq:dim_sum}
        \operatorname{dim}\left(\mathscr{W}^{\perp_{\operatorname{s}}}\right)+\operatorname{dim}\left(\mathscr{W}\right)= \operatorname{dim}\left(\mathscr{V}\right).
    \end{align}
    See \cite[Proposition~1.13]{degosson}.
    A linear subspace $\mathscr{W}$ of $\mathscr{V}$ is said to be a \emph{symplectic subspace} if the intersection of $\mathscr{W}$ and $\mathscr{W}^{\perp_{\operatorname{s}}}$ is the zero subspace, or  equivalently, $\omega$ restricted to $\mathscr{W}$ is also non-degenerate.
    A subspace $\mathscr{W}' \subseteq \mathscr{V}$ is said to symplectically orthogonal to $\mathscr{W}$ if $\mathscr{W}' \subseteq \mathscr{W}^{\perp_{\operatorname{s}}}$.

    % The following relationship between the Poisson bracket of two quadratic forms $Q, R$ and their Hamilton maps will be useful
    % \cite[Lemma~2]{pravda2008contraction}\footnote{It is Lemma~3.2 in the arXiv version of the paper. Also, there is a factor of $2$ that appears in \cite{pravda2008contraction} but not in our case. This is because we consider the factor $\frac{1}{2}$ in our definition of the polar identity associated with a quadratic form.}: 
    % \begin{align}\label{eq:poisson_brac_hamilton}
    %     H_{\{Q, R\}} = -[H_Q, H_R],
    % \end{align}
    % where $[H_Q, H_R]\coloneqq H_QH_R-H_RH_Q$ is the commutator of $H_Q$ and $H_R$.
    
    We now discuss interpretations of symplectic orthogonal projection and some of the results for quadratic forms on general symplectic spaces.
    We emphasize that the translations of the results for a quadratic form $Q$ are obtained by the corresponding symmetric matrix of $\Phi_Q$ in a symplectic basis.
 \subsection{General Williamson's theorem (Theorem~\ref{thm:williamson_general})}
    Let $Q$ be a quadratic form on a symplectic space $(\mathscr{V}, \omega)$.
    Theorem~\ref{thm:williamson_general} states that there exists a symplectic basis 
    $\{p_1,\ldots, p_n, q_1, \ldots, q_n\}$ of $\mathscr{V}$, and real numbers $\mu_1,\ldots, \mu_n$ such that for all $(x_1,\ldots, x_n, y_1,\ldots, y_n) \in \mathds{R}^{2n}$,
    \begin{align} \label{eq:general_will_quadratic_form}
        Q\left(\sum_{i=1}^n (x_i p_i + y_i q_i)\right)
            = \sum_{i=1}^n \mu_i \left(x_i^2 + y_i^2 \right),
    \end{align}
    if and only if there exist symplectic subspaces $\mathscr{W}_-$, $\mathscr{W}_0$, $\mathscr{W}_+$ of $\mathscr{V}$ with dimensions $\nu(Q), \xi(Q), \pi(Q)$, respectively, such that
    \begin{itemize}
            \item [$\circ$] $\mathscr{W}_-$, $\mathscr{W}_0$, $\mathscr{W}_+$ are pairwise symplectically orthogonal to each other,
            \item [$\circ$] these subspaces are invariant under the Hamiltonian operator $H_Q$, and
            \item [$\circ$] $Q$ takes strictly negative values on $\mathscr{W}_-$, vanishes on $\mathscr{W}_0$, and it takes strictly positive values on $\mathscr{W}_+$.
    \end{itemize}
    Furthermore, the numbers $\mu_1,\ldots, \mu_n$ are unique.
    Moreover, $\pm \mu_1,\ldots, \pm \mu_n$ are the eigenvalues of $i H_Q$ over the complexification of $\mathscr{V}$\footnote{The complexification of $\mathscr{V}$ is a complex vector space $\mathscr{V}_{\mathds{C}} \coloneqq \mathscr{V} \oplus i  \mathscr{V}$ with the vector addition and scalar multiplication defined in a natural way. That is, for $u_1, u_2, v_1, v_2 \in \mathscr{V}$ and $\alpha, \beta \in \mathds{R}$
    \begin{align}
        (u_1+ i  v_1) + (u_2+i  v_2) &\coloneqq (u_1+u_2)+i  (v_1+v_2), \\
        (\alpha + i  \beta)(u_1+ i  v_1) &\coloneqq (\alpha u_1-\beta v_1)+i  (\beta u_1+\alpha v_1).
    \end{align}
    Every real linear map $H: \mathscr{V} \to \mathscr{V}$ can be extended to a complex linear map $H:\mathscr{V}_{\mathds{C}} \to \mathscr{V}_{\mathds{C}}$ as
    \begin{align}
        H(u+i  v) \coloneqq H(u)+i  H(v), \qquad u,v \in \mathscr{V}.
    \end{align}
    }.
 \subsection{Description of symplectic eigenvalues for $\operatorname{EigSpSm}(2n)$ (Theorem~\ref{thm:williamson_gen_special})}
    Let $Q$ be a quadratic form on a symplectic space $(\mathscr{V}, \omega)$, and let
    $J$ be a complex structure on $\mathscr{V}$ compatible with $\omega$.
    Let $A$ be the symmetric matrix representing the bilinear form $\Phi_Q$ in a $J$-orthosymplectic basis, and suppose that $A$ belongs to $\operatorname{EigSpSm}(2n)$\footnote{This property is independent of the choice of the $J$-orthosymplectic basis. This is because an automorphism taking a $J$-orthosymplectic bases to another $J$-orthosymplectic basis is given by an orthosymplectic matrix.}.
    It is straightforward to see that $Q$ can be brought into Williamson's normal form \eqref{eq:general_will_quadratic_form} in a symplectic basis and the symplectic eigenvalues of $A$ are $\mu_1,\ldots, \mu_n$.
    The conclusion of Theorem~\ref{thm:williamson_gen_special} holds for $A$.

\subsection{Symplectic orthogonal projection}
    A symplectic orthogonal projection in a general symplectic space $\left(\mathscr{V}, \omega \right)$ is a projection or \emph{idempotent map} $\Pi: \mathscr{V} \to \mathscr{V}$ such that
    \begin{itemize}
        \item $\ker(\Pi)$ is a symplectic subspace of $\mathscr{V}$, and
        \item $\operatorname{range}(\Pi)= \ker(\Pi)^{\perp_{\operatorname{s}}}$.
    \end{itemize}
    The statement of Proposition~\ref{prop:transpose-symp-subspace-projection} holds with the adjoint operator of $\Pi$ with respect to the inner product $g_J$ on $\mathscr{V}$ induced by a complex structure $J$ compatible with $\omega$.
    Indeed, suppose $J$ is a complex structure on $\mathscr{V}$ compatible with $\omega$.
    Let $\Pi$ be a symplectic orthogonal projection, and let $\Pi^{\perp_{\operatorname{s}}}$ be the adjoint of $\Pi$ with respect to the inner product $g_J$. 
    Since $\Pi$ is idempotent, $\Pi^{\perp_{\operatorname{s}}}$ is also idempotent.
    We have for arbitrary $u \in \ker(\Pi)$ and $v \in \mathscr{V}$ that
    \begin{align}
        g_J(\Pi^{\perp_{\operatorname{s}}}(Ju), v)
            &= g_J(Ju, \Pi(v)) \\
            &= \omega(Ju, J\Pi (v)) \\
            &= \omega(u, \Pi (v)) \\
            &=0.
    \end{align}
    This implies that $\Pi^{\perp_{\operatorname{s}}}(Ju)=0$ and hence $J(\ker(\Pi)) \subseteq \ker(\Pi^{\perp_{\operatorname{s}}})$.
    Since $J$ is an automorphism and $\rank(\Pi^{\perp_{\operatorname{s}}})=\rank(\Pi)$, we thus conclude that 
    \begin{align}
        \ker(\Pi^{\perp_{\operatorname{s}}})=J(\ker(\Pi)).
    \end{align}
    Also, 
    \begin{align}
        \omega(Ju, \Pi^{\perp_{\operatorname{s}}}(v))
            &= - \omega(u, J(\Pi^{\perp_{\operatorname{s}}}(v))) \\
            &= -g_J(u, \Pi^{\perp_{\operatorname{s}}}(v)) \\
            &= -g_J(\Pi(u), v) \\
            &=0.
    \end{align}
    This implies that 
    \begin{align}
        \operatorname{range}(\Pi^{\perp_{\operatorname{s}}}) 
            &\subseteq J(\ker(\Pi))^{\perp_{\operatorname{s}}} \label{eq:symp_perp_range} \\
            &=(J(\ker(\Pi))^{\perp_{\operatorname{s}}} \\
            &= \ker(\Pi^{\perp_{\operatorname{s}}})^{\perp_{\operatorname{s}}}=J(\operatorname{range}(\Pi)).
    \end{align}
    The rank-nullity theorem, combined with the relation \eqref{eq:dim_sum}, implies that the inclusion in \eqref{eq:symp_perp_range} cannot be proper.
    We have thus proved that $\Pi^{\perp_{\operatorname{s}}}$ is a symplectic orthogonal projection whose range is given by $J(\operatorname{range}(\Pi))$.

%%%%%%%%%%%%%%%%%%%%%%%%%%%%%%%%%%%%%%%%%%%%%%%%%%%%%%%%%%%%%%%%%%

\section*{Acknowledgments}
    The author thanks the anonymous referee for their suggestions for improving the readability of the paper and for highlighting the geometrical aspects of some of the results.
     The author acknowledges supports from the NSF under grant no.~2304816, AFRL under agreement no.~FA8750-23-2-0031, and FRS Project No.~MISC~0147. 
     The author thanks Prof.~Tanvi Jain for insightful discussions.

\begin{backmatter}

\bibliographystyle{unsrtnat}
\bibliography{RefBT}
\end{backmatter}
\printaddress

\end{document}